\documentclass[a4]{amsart} \usepackage{verbatim} 
\usepackage{amssymb}
\usepackage{amscd}

\newcommand{\halfline}{{\RR_{\geq 0}}}
\newcommand{\so}{P}
\newcommand{\nbt}{{\textbf{B}}}
\newcommand{\nft}{{\textbf{F}}}
\newcommand\half{\frac1{2}}
\newcommand\Dist{C^{-\infty}}

\newcommand{\FT}{{\mathcal{F}}}

\newcommand{\tA}{{\tilde{A}}}
\newcommand{\cH}{{\mathcal{H}}}
\renewcommand\MR[1]{}

\theoremstyle{plain}
\newtheorem{theorem}{Theorem}[section]
\newtheorem{proposition}[theorem]{Proposition}
\newtheorem{lemma}[theorem]{Lemma}

\theoremstyle{definition}

\theoremstyle{remark}
\newtheorem*{remark}{Remark}

\theoremstyle{example}

\newlength{\KDlen}
\newlength{\sclen}
\newlength{\tlen}
\newlength{\slen}
\newlength{\qsclen}

\newlength{\omegalen}

\numberwithin{equation}{section}

\DeclareMathSymbol{\leqslant}{\mathrel}{AMSa}{"36}
\DeclareMathSymbol{\geqslant}{\mathrel}{AMSa}{"3E}
\DeclareMathSymbol{\gtrless}{\mathrel}{AMSa}{"3F}
\DeclareMathSymbol{\lessgtr}{\mathrel}{AMSa}{"37}

\renewcommand{\leq}{\leqslant}
\renewcommand{\geq}{\geqslant}

\def\sch/{{Schr\"odinger}} 

\def\psido/{{$\Psi$DO}}

\newcommand{\SC}{\ensuremath{\mathrm{sc}}}
\newcommand{\QSC}{\ensuremath{\mathrm{qsc}}}

\newcommand{\tf}{\text{sf}}

\newcommand\restrictedto{\!\!\upharpoonright}

\newcommand{\RR}{{\mathbb{R}}}

\newcommand{\NN}{{\mathbb{N}}}

\newcommand{\loc}{{\text{loc}}}

\newcommand\CI{{\mathcal{C}}^{\infty}}

\newcommand{\CdotI}{{\dot{\mathcal{C}}^\infty}}

\newcommand{\h}{{\frac{1}{2}}}
\newcommand{\ep}{{\epsilon}}

\DeclareMathOperator{\wfront}{{WF}}

\newcommand{\WF}[1][\mbox{}]{{\wfront^{#1}}}
\newcommand{\WFsc}[1][\mbox{}]{{\wfront_{\SC}^{#1}}}

\newcommand{\opWF}{{\wfront'}}

\newcommand{\Tsc}[1][\mbox{}]{{\settowidth{\tlen}{$T$}
    \settowidth{\sclen}{{\scriptsize \SC}} \hspace{\sclen}
    T_{#1}^{\hspace{-\sclen}\hspace{-\tlen}
      \SC\hspace{\tlen}}}}

\newcommand{\Tscstar}[1][\mbox{}]{{\settowidth{\tlen}{$T$}
    \settowidth{\sclen}{{\scriptsize
        \SC}}
 \hspace{\sclen}T_{#1}^{\hspace{-\sclen}\hspace{-\tlen}
      \SC
 \hspace{\tlen}
      \ast}}}

\newcommand{\Tqsc}{{\settowidth{\tlen}{T}
    \settowidth{\qsclen}{{\scriptsize
 \QSC}} \hspace{\qsclen}
    T^{\hspace{-\qsclen}\hspace{-\tlen}
      \QSC \hspace{\tlen}}}}

\newcommand{\Tqscstar}[1][\mbox{}]{{\settowidth{\tlen}{T}
    \settowidth{\qsclen}{{\scriptsize \QSC}} \hspace{\qsclen}
    T_{#1}^{\hspace{-\qsclen}\hspace{-\tlen} \QSC
 \hspace{\tlen}
      \ast}}}

\newcommand{\nqsc}[1][\cdot]{N_{\QSC}}

\DeclareMathOperator{\supp}{{supp}}

\newcommand{\pa}{{\partial}}

\newcommand{\abs}[1]{{\left\lvert{#1}\right\rvert}}
\newcommand{\norm}[1]{{\left\lVert{#1}\right\rVert}}

\newcommand{\Psisc}{{}^{\SC}\Psi}

\newcommand{\Omegasc}{{\settowidth{\omegalen}{$\Omega$}
    \settowidth{\sclen}{{\scriptsize \SC}} \hspace{1.1\sclen}\Omega^{\hspace{-1.1\sclen}\hspace{-\omegalen}
      \SC}\hspace{\omegalen}}}

\newcommand{\ocal}{{\mathcal{O}}}

\newcommand{\Hsc}[1][\mbox{}]{H_{\SC}^{#1}}

\newcommand\Lt{\tilde
L}
\newcommand\mf{\operatorname{mf}}
\renewcommand\tf{\operatorname{sf}}

\newcommand\XtS{X_{\operatorname{tS}}}
\newcommand\cU{\mathcal{U}}
\newcommand\cUZ{\mathcal{U}_Z}
\newcommand\cUt{{\tilde{\mathcal{U}}}}
\newcommand\cRt{{\tilde{\mathcal{R}}}}
\newcommand\cE{\mathcal{E}}
\newcommand\cR{\mathcal{R}}
\newcommand\cG{\mathcal{G}}
\newcommand\XXt{X^2_t}
\newcommand\Vphi{
    {}^\phi \mathcal{V}}
\newcommand\TphiXXt{ {}^\phi T X^2_t
  }
\newcommand\TstarphiXXt{ {}^\phi T^* X^2_t}
\newcommand\Lphi{
    {}^\phi L}
\newcommand\HDsc{ {}^{\operatorname{sc}} \Omega^{1/2}}
  \newcommand\ie{i.e.}
 \newcommand\scH{{}^{\operatorname{sc}}H}

\newcommand\contact{S}

\newcommand\sub{\operatorname{sub}}
\newcommand\Fourier{\mathcal{F}}
\newcommand\Openset{\textbf{O}}

\title[Schr\"odinger propagator]{The Schr\"odinger propagator for scattering metrics}
\author{Andrew Hassell}
\address{Department of Mathematics, Australian National University, Canberra, 0200 ACT Australia}
\author{Jared Wunsch}
\address{Department of Mathematics, Northwestern University, Evanston IL USA}
\thanks{This research was supported in part by a Fellowship and a Linkage grant from the Australian Research Council (A. H.) and by NSF grant DMS-0100501 (J. W.)}

\begin{document}
\maketitle


\begin{abstract} Let $g$ be a scattering metric on a compact
  manifold $X$ with boundary, i.e., a smooth metric giving the interior $X^\circ$ the structure of a complete Riemannian manifold with  asymptotically conic ends. An example is any compactly supported perturbation of the standard metric on $\RR^n$. 
  Consider the operator $H = \half \Delta + V$, where $\Delta$ is the
  positive Laplacian with respect to $g$ and $V$ is a smooth
  real-valued function on $X$ vanishing to second order at $\pa X$.  Assuming that $g$ is non-trapping, we
  construct a global parametrix $\mathcal{U}(z, w,t)$ for the kernel of the Schr\"odinger
  propagator $U(t) = e^{-i t H}$, where $z, w  \in X^{\circ}$  and $t \neq 0$. The parametrix is such that the difference between $\mathcal{U}$ and $U$ is smooth and
  rapidly decreasing both  as $t \to 0$ and as $z \to \pa X$, uniformly for $w$ on compact subsets of $X^{\circ}$. Let $r =
  x^{-1}$, where $x$ is a boundary defining function for $X$, be an
  asymptotic radial variable, and let $W(t)$ be the kernel
  $e^{-ir^2/2t}U(t)$. Using the parametrix, we show that $W(t)$ belongs to a class of
  `Legendre distributions' on $X \times X^{\circ} \times
  \halfline$ previously considered by Hassell-Vasy.
When the metric is trapping, then the parametrix construction goes through microlocally in the non-trapping part of the phase space. 

  We apply this result to obtain a microlocal characterization of the singularities of $U(t) f$, for any tempered distribution $f$ and any fixed $t \neq 0$, in terms of the oscillation of $f$ near $\pa X$. 
If the metric is non-trapping, then we obtain a complete characterization; more generally we need to assume 
that $f$ is microsupported in the nontrapping part of the phase space. This generalizes results of Craig-Kappeler-Strauss and Wunsch. 
\end{abstract}
 
\section{Introduction}\label{section:intro}
Let $(X,g)$ be a scattering manifold of dimension $n$. Thus, $X$ is a compact $n$-dimensional manifold with 
boundary, and $g$ is a metric in the interior of $X$ taking the form
$$
g = \frac{dx^2}{x^4} + \frac{h}{x^2}
$$ near the boundary. Here $x$ is a boundary defining function, $y$ are
local coordinates on the boundary $Y = \pa X$ extended to $X$ and $h$ is a
2-cotensor that restricts to a metric (i.e., is positive definite) on $Y$. We
shall assume that $x$ is a globally defined, smooth function on $X$ that
vanishes only at $\pa X$, and denote $r = x^{-1}$, which is analogous to
the radial variable on Euclidean space. The metric takes a more familiar
form when written in terms of $r$:
$$
g = dr^2 + r^2 h(r^{-1}, y, dr, dy),
$$ 
where $h$ is smooth in the first variable and $y$ (smoothness in $x =1/r
\in [0, \epsilon)$ is of course a much stronger condition than smoothness
in $r \in (\epsilon^{-1}, \infty)$). In fact, one can choose coordinates
locally so that $g$ takes the form
\begin{equation}
g = dr^2 + r^2 h(r^{-1}, y, dy) = \frac {dx^2}{x^4} + \frac{h(x,y,dy)}{x^2};
\label{better}\end{equation}
here $h$ is a metric in the $y$ variables depending parametrically (and
smoothly) on $x$ (see \cite{Joshi-SaBarreto2} for a proof). Thus, $g$ is {\it asymptotically
conic}; it approaches the conic metric $dr^2 + r^2 h_0$, where $h_0 = h(0,
y, dy)$, as $r \to \infty$. The boundary $\pa X$ is then at geometric infinity, with each point of $\pa X$ representing an asymptotic direction of geodesics. Euclidean space, with its standard metric or any compactly supported perturbation of it, is an example. 
We shall assume henceforth that coordinates have been
chosen so that the representation \eqref{better} for the metric holds.
Let $\Delta$ denote 
the positive Laplacian with respect to $g$, and let $H = \half \Delta + V$, 
where $V\in x^2\CI(X;\RR).$ Thus $V$ is a short-range potential\footnote{The case of gravitational long-range potentials in the sense of \cite{MR2002i:58037} can also be treated.}. We consider the time-dependent Schr\"odinger equation
\begin{equation}
(D_t + H ) u(z, t)=0, \quad u(z, 0) = f(z) \in \Dist(X), 
\label{sch-eqn}\end{equation}
where $z \in X$, $t \in \RR$, $D_t = -i \frac{\pa}{\pa t}$ and $f$ is a given distribution on $X$.  Let
$U(t) = e^{-i t H}$ be the propagator for $H$. We wish to construct a
parametrix $\mathcal{U}(t)$ for $U(t)$ which captures all the singularities
of $U$; in particular, we want $\mathcal{U}(t) - U(t)$ to have a kernel
which vanishes rapidly both as $t \to 0$ and as we approach $\pa X$. However, to simplify the
construction we shall let only one of the variables $(z,w) \in X^2$, say
the left variable $z$, approach infinity, while we shall only require uniformity over compact subsets in the $w$ variable. 

To state our main theorem, we choose a function $\phi \in \CI(X)$
which is zero in a neighbourhood of the boundary of $X$, we let $r =
x^{-1}$ as above, and we define the kernel $W(t) = e^{-ir^2/2t}
U(t)$. We remark that we regard $U(t)$ and $W(t)$ as acting on
half-densities, so the kernels contain a Riemannian half-density
factor $|dg_z dg_w dt|^{1/2}$ in each variable. Then our main result is

\begin{theorem}\label{main} Assume that the metric $g$ is non-trapping. Then the kernel $\phi(w) W(t)$ is a fibred-scattering Legendrian distribution on
$X \times X^\circ \times \RR_+$ of order $(\frac{3}{4}, \frac{1}{4})$ in the sense of
\cite{MR2001d:58034}, associated to the Legendre submanifold $\Lphi$
defined in Lemma~\ref{lemma:fibredLeg}.
\end{theorem}

\begin{remark} In the trapping case, one can still construct a parametrix microlocally in the non-trapping region; see Section~\ref{section:trapping}. 
\end{remark}

In more prosaic language, Theorem~\ref{main} says that the propagator is given by oscillatory integrals
of certain rigidly prescribed forms. Near $t = 0$, and with both variables
$w,z$ away from the boundary of $X$, these take the form of a Legendrian
distribution (see Melrose-Zworski \cite{MR96k:58230})
$$
t^{-n/2 - k/2} (2\pi)^{-k} \int_K e^{i\Phi(z,w,v)/t} a(z,w,v,t) dv,
$$
where the integral is over a compact set $K \subset \RR^k$. Here, $\Phi$ and $a$ are smooth in all their variables. Near the diagonal, $\Phi$ is given by
$\operatorname{dist}(w,z)^2/2$, as in the flat Euclidean case, no $v$ variables are required, and there is no integral. Associated to $\Phi$ is a Legendrian submanifold of $T^*X \times T^* X \times \RR$, here given explicitly by
\begin{equation*}
L = \{ (w, \xi, z, \zeta, \tau) \mid (z, \zeta) = \exp_{g}(w, \xi), \ \tau = \frac{\abs{\xi}^2 }{2}
\}
\end{equation*}
(this is a somewhat non-invariant description, for precise details see
\S\ref{section:flowout}).  The function $\Phi$ becomes non-smooth
outside the injectivity radius; geometrically this corresponds to the
Legendrian $L$ becoming non-projectable (i.e.\ the projection from the
Legendrian to the base $X \times X$ is no longer a
diffeomorphism). The Legendrian, however, remains perfectly
smooth; it is no longer parametrized by a function just of $w$ and
$z$, but one needs extra variables $v$ (precisely, one needs at least
$k$ extra variables locally if the kernel of the projection from the
Legendrian to the base has dimension $k$).

When $z$ approaches the boundary of $X$, we use local coordinates $(x,y)$
as described above, and sometimes use $r = x^{-1}$, the asymptotic `radial'
variable. In this region the propagator takes a more complicated form. We show that $W(t)$
is a finite sum of oscillatory integrals of the form
\begin{equation}
t^{-n/2 - k/2} (2\pi)^{-k} \int_K e^{i\psi(x,y,w,v)/xt} a(x,y,w,v,t) dv,
\end{equation}
where again $\psi$ and $a$ are smooth in all their variables, and one integrates over a compact set $K \subset \RR^k$. Notice that $\psi$ depends on all variables apart from $t$. 

\

Next, we address the question of determining the wavefront set of $u(\cdot,
t)$ at a fixed nonzero (say, positive) time $t$ in terms of the wavefront
set of the initial data $f$.  We are particularly interested in the case of
`interior singularities', lying above a point in the interior of $X$. It is
well known that equation \eqref{sch-eqn} has infinite speed of propagation, so to
answer this question we must look for singularities of $f$ at the boundary of $X$. Dually, since $U(t) = U(-t)^*$, we can consider interior singularities of $f$ and determine the singularities of $u(t)$ that they produce.  It is known
from \cite{Wunsch1} that if there is an interior singularity $(w, \eta)$ for
$f$, then for all positive times $t$ there are quadratic oscillations of
frequency $1/2t$ in $u(\cdot, t)$ at infinity in the asymptotic direction
of the geodesic emanating from $(w, \eta)$. Roughly speaking, the quadratic
oscillations look like $e^{ir^2/2t}$ in a conic neighbourhood of the
asymptotic direction $y$; more precisely, this result is phrased in terms
of the `quadratic scattering wavefront set' (see \cite{Wunsch1}). 

The limitation of the result of \cite{Wunsch1} is that different singularities along a single geodesic, or along a different geodesic with the same asymptotic direction, produce identical quadratic oscillation, hence consideration of the quadratic wavefront set alone will not result in precise propagation results. To analyse the finer structure of $u(\cdot, t)$, we divide   $u(\cdot, t)$ by the explicit quadratic oscillatory factor $e^{ir^2/2t}$, and
find that the resulting function has oscillations which are {\it linear} in $r$ that contain the
desired information on the location of the interior singularity. The presence of linear oscillations of a function $f$ is measured by the \emph{scattering wavefront set} $\WFsc(f)$ (see \cite{MR95k:58168}), a closed subset of $\Tscstar[\pa X] X$. (This bundle, defined in \S\ref{section:definitions}, is  the restriction to $\pa X$ of an $n$-dimensional bundle over $X$ which is a compressed and scaled version of the cotangent bundle $T^*X$).  We show that the asymptotics of geodesic flow on $X$ determine two
contact transformations $S_f$ and $S_b$, which we call the forward and backward sojourn relations, from
the sets $\nft$, resp. $\nbt \subset S^* X^{\circ}$ consisting of points in the cosphere bundle which are not forward, resp. backward trapped under geodesic flow, to $\Tscstar[\pa X] X$. They are related by $\nft = -\nbt$ and $S_f(\zeta) =  -S_b(-\zeta)$. The 
definition of $S_f$ in local coordinates is as follows:
Let $\gamma(s)$ be the arclength-parametrized geodesic
emanating from $(w, \eta) \in S^* X^{\circ}$, let $y = (y_1, \dots, y_{n-1})$ be local coordinates on $\pa X$ and let $(y, \nu, \mu = (\mu_1, \dots, \mu_{n-1}))$ be the induced coordinates on $\Tscstar[\pa X] X$ as in \eqref{numu}.  Then $S_f(w, \eta) = (y_0, \nu, \mu)$ iff 
$$
y_0=\lim_{s\to +\infty} \gamma(s) \in \pa X
$$
is the asymptotic direction of the geodesic,
\begin{equation}
\nu = \lim_{s \to +\infty} \big(s - r(\gamma(s)) \big), \quad r = \frac1{x}
\label{soj-defn}\end{equation}
 is the `sojourn time' of the geodesic, and 
\begin{equation}
\mu_i = \sum_j h_{ij} \lim_{s \to +\infty}
(y_0 - y(\gamma(s)))_j)/x
\label{pencil-coord}\end{equation}
measures its angle of approach to $\pa X.$
(That the point $y_0\in \pa X$ should be thought of as the asymptotic
\emph{direction} of the geodesic can be seen easily in the special
case when $X$ is the radial compactification of $\RR^n.$) Note that
the sojourn time, thought of as depending on a geodesic and a point
along it, is closely related to the classical Busemann function of
differential geometry; see for example \cite{Schoen-Yau}, chapter 1, section 2. We use the term `sojourn relation' for $S_f$ since the coordinate $\nu$ is analogous to the sojourn time considered by Guillemin \cite{Gu}. 

\begin{theorem}\label{WF-thm} Suppose $f \in \Dist(X)$, and let $u(\cdot, t) = e^{-itH} f$. Let $\zeta=(w, \eta) \in \nbt \subset S^*X^\circ$ be an non-backward-trapped point of $S^*X^\circ$.
Then for any fixed $t > 0$,
\begin{equation}
\zeta \in \WF(u(\cdot, t))
\text{  if and only if } \ 
\frac1{t} S_b(\zeta) \in \WFsc(e^{ir^2/2t} f).
\label{sc-wf}\end{equation}
Here the factor $t^{-1}$ acts by scaling the fibre variables. 
Similarly, if $\zeta \in \nft$, 
then for  $t < 0$, 
\begin{equation}
\zeta \in \WF(u(\cdot, t))
\text{  if and only if } \ 
\frac1{|t|} S_f(\zeta) \in \WFsc(e^{ir^2/2t} f).
\label{sc-wf-2}\end{equation}
\label{theorem:propagation}
\end{theorem}

\begin{remark} The condition on the right hand side of  \eqref{sc-wf} or \eqref{sc-wf-2} is not manifestly coordinate invariant. The sojourn relation changes under a change of coordinates, but so does the scattering wavefront set of $e^{ir^2/2t} f$, in such a way that the condition as a whole is coordinate invariant. A manifestly invariant description may be given in terms of the affine bundle of Lemma~\ref{Lt}. 
\end{remark}

\noindent \emph{Example.} Consider the free Hamiltonian $H = \frac1{2} \Delta$ on $\RR^n$. The propagator in this case is 
$$
U(t) = (2\pi t)^{-n/2} e^{i|z-w|^2/2t}, \quad \text{whence} \quad 
W(t) = (2\pi t)^{-n/2} e^{-iz \cdot w/t} e^{i|w|^2/2t}.
$$
It is not hard to check that $W(t)$ satisfies the conclusion of Theorem~\ref{main} in this case. Consider initial data $$f = (2\pi T)^{-n/2} e^{-i|z-w|^2/2T}, \quad T > 0.$$ Then the solution at time $t = T$ is a delta function centred at $w$, hence its wavefront set is $$\{ (w, \hat z) \mid \hat z \in S^{n-1} \}.$$ In the free case, the backward sojourn relation is given by $$S_b(w, \hat z) = \big(y = -\hat z, \, \nu = - w \cdot \hat z, \, \mu = w - (w \cdot \hat z) \hat z \big),$$ and we easily check that \eqref{sc-wf} holds in this case. 

\vskip 5pt

Most previous work on Schr\"odinger parametrices has focused on the case of
flat space with a potential perturbation, where the geometric situation is
substantially simpler.  Very detailed parametrix constructions in this
setting have been made by Fujiwara \cite{Fujiwara1}, Zelditch
\cite{Zelditch2}, Tr\`eves \cite{Treves1} and Yajima \cite{Yajima1}.  In the
case of curved space, very little was known.  Kapitanski-Safarov 
\cite{KapSaf1} have shown that on $\RR^n$ with a compactly-supported,
nontrapping potential perturbation, the fundamental solution is smooth for
$t>0,$ and have exhibited a parametrix modulo $\CI(\RR^n)$ \cite{KapSaf2}.
Such a parametrix, however, is not sufficiently specified at infinity to
yield results about smoothness of the solution of the general Cauchy
problem at $t>0.$

Regularity results for the Schr\"odinger propagator for non-flat metrics in the form of Strichartz estimates have been obtained recently by Staffilani-Tataru \cite{Staff-T} and Burq-G\'erard-Tzvetkov \cite{BGT}. 
Staffilani and Tataru \cite{Staff-T} proved Strichartz estimates for $e^{-itH} f$ where $H$ is the Laplacian of a $C^2$, compactly supported, nontrapping perturbation of the standard metric on $\RR^n$, using the FBI transform and Littlewood-Paley decompositions to handle the rough metric. Burq, G\'erard and Tzvetkov obtained Strichartz estimates, with a loss of derivatives compared to the flat Euclidean case, for compact manifolds or perturbations of the Laplacian on $\RR^n$, without any nontrapping assumption.

Various authors have considered the question of determining the singularities of $u = e^{-itH} f$ in terms of $f$. 
The first results about microlocal smoothness of $u(t)$ for $t>0$ and general initial data
were those of Craig-Kappeler-Strauss \cite{CKS}, who showed that, on
an asymptotically Euclidean space, decay of the initial data in a
microlocal incoming cone yields microlocal smoothness along the whole
pencil of geodesics emanating from that cone for all $t>0.$ This
result was refined in \cite{Wunsch1}, where the second author showed
that absence of the quadratic-scattering wavefront set (see 
Section~\ref{qsc-str}) allowed one to conclude
absence at varying \emph{times} and along varying pencils of
geodesics.  This approach, while it specified in terms of the Cauchy
data \emph{when} and \emph{in what direction} singularities might
appear in $X^\circ,$ failed to say anything about \emph{where} they
might land.  These propagation results have been extended to
the analytic category by Robbiano-Zuily \cite{R-Z}. 

Our Theorem~\ref{WF-thm} gives a complete solution to the propagation problem in the case when the metric $g$ is nontrapping, and in general, a complete characterization of the singularities of $u(t)$ in the non-backward-trapped set $\nft$ for $t>0$, and in the non-forward-trapped set $\nbt$ for $t < 0$. Our results imply those of \cite{Wunsch1}, and hence those of \cite{CKS}\footnote{although we require more decay of our potential than is assumed in \cite{CKS}.}, since the hypothesis on $f$ required in \cite{Wunsch1} for
microlocal smoothness of $u(t)$ along all geodesics emanating
from $y \in \pa X$ implies that in fact $(y,\nu,\mu) \notin \WFsc
e^{ir^2/2 t} f$ for \emph{all} $(\nu,\mu)$. Hence
Theorem~\ref{theorem:propagation} yields the main boundary to interior
propagation result of \cite{Wunsch1} as a special case.

We thank Richard Melrose, Andr\'as Vasy and Steve Zelditch for useful conversations, and the Erwin Schr\"odinger Institute, the Mathematics Department at SUNY Stonybrook and the Mathematical Sciences Institute at the Australian National University for their hospitality. 


\section{Contact structures and Legendrian distributions}\label{section:definitions}
We recall the definition of various structures associated to manifolds
with boundary and corners needed in this paper. For further details,
see Melrose \cite{MR95k:58168}, Melrose-Zworski \cite{MR96k:58230},
Hassell-Vasy \cite{MR2002i:58037,MR2001d:58034}.

\subsection{Scattering structure and Legendrian
  distributions}\label{subsec:Leg}

Let $X$ be a $d$-dimensional manifold with boundary, and let $x$ be a
boundary defining function for $X$.  We identify a collar neighborhood
of $\pa X$ with $\pa X \times [0,\ep)$ so that specifying $y_1,\dots
  y_{d-1}$ local coordinates in $\pa X$ gives local coordinates
  $(x,y)$ for $X$. The Lie algebra of scattering vector fields
  $\mathcal{V}_{\SC}(X)$ consists of vector fields of the form $a x^2
  \pa_x + \sum x b_i \pa_{y_i}$ with $a, b_i \in \CI(X)$. Such
  vector fields can be described as the set of $\CI$ sections of a
  vector bundle $\Tsc X$. The dual of this bundle we denote $\Tscstar
  X$; sections of it are locally spanned over $\CI(X)$ by $dx/x^2$ and
  $dy_i/x$.  Hence any point $q \in \Tscstar[x,y] X$ has a unique
  expression
\begin{equation}
q = \nu d(\frac1{x}) + \sum \mu_i \frac{dy_i}{x}
\label{numu}\end{equation}
which yields local coordinates $(x, y, \nu, \mu)$ for $\Tscstar X$,
$\nu$ and $\mu$ being linear on each fibre. We say that a half-density
$\alpha$ on $X$ is a scattering half-density if it is a smooth and
non-vanishing section of the bundle $(\wedge^n (\Tscstar X))^{1/2}$;
such sections have the form $a |dx dy_1 \dots dy_{d-1}/x^{d+1}|^{1/2}$,
where $a \in \CI(X)$ is nonvanishing. The restriction of $\Tscstar X$ to the boundary of $X$ is denoted $\Tscstar[\pa X] X$. 

We define $\CdotI(X) = \bigcap_{l \geq 0} x^l C^\infty(X)$, with its natural Fr\'echet topology, and denote by $\Dist(X)$ its topological dual. We sometimes refer to these as the space of Schwartz functions and the space of tempered distributions on $X$, by analogy with Euclidean space.

Now we recall some facts about the scattering calculus
$\Psisc^{m,l}(X; \Omegasc^{1/2}(X))$ acting on half-densities, which
is indexed by two orders $(m,l)$, the interior order $m$ (which for a
differential operator is the order of the highest derivative that
occurs) and the boundary order $l$. The half-density factor will be
understood from now on, and dropped from notation. The space $\Psi^{m,l}(X)$ is the same as $x^l \Psi^{m,0}(X)$, and $\Psi^{m,0}(X)$ is a `microlocalization' of the scattering differential operators of order $m$ on $X$; it contains in particular all $m$th order differential operators generated over $\CI(X)$ by $\mathcal{V}_{\SC}(X)$. Operators $P \in
\Psisc^{m,l}(X)$ are determined up to $\Psisc^{m, l+1}(X)$ by the
boundary symbol $p$, which is a smooth function on the boundary of the
scattering cotangent bundle, $\Tscstar[\pa X] X$, and up to
$\Psisc^{m, l+2}(X)$ by the boundary symbol $p$ together with the
boundary subprincipal symbol, which is again a smooth function
$p_{\sub}$ on $\Tscstar[\pa X] X$. The operator $P$ is said to be elliptic at $q \in \Tscstar[\pa X] X$ if $p(q) \neq 0$. The scattering wavefront set
$\WFsc(u)$ of a distributional half-density $u$ is defined by the
condition that $q \in \Tscstar[\pa X] X$ is \emph{not} in $\WFsc(u)$ iff there is an $A \in
\Psisc^{0,0}(X)$ such that $A$ is elliptic at $q$ and $Au \in
\CdotI(X),$ the space of Schwartz functions on $X^\circ.$ The
scattering wavefront set of $u$ is always a closed subset of
$\Tscstar[\pa X] X$.  We also have a scale of Sobolev spaces
$\Hsc[m](X)$, defined by $u \in \Hsc[m](X)$ iff $V_1 \cdots V_m u \in L^2(X)$
for all $V_1, \dots , V_m \in \mathcal{V}_{\SC}(X)$, or equivalently, if $Pu \in L^2(X)$ for all $P \in \Psi^{m,0}(X)$. (Here $L^2(X)$ is defined with respect to the Riemannian measure $dg$.)

For the purposes of this paper, we will often take the manifold with
boundary to be $X^\circ \times X^\circ \times \halfline$ (this space is not compact,
but that is irrelevant here). In that case the boundary defining
function is $t$ and local coordinates on the boundary will be denoted
$(z,w)$, where $z \in \RR^n$ is a local coordinate for the first
factor and $w \in \RR^n$ is a local coordinate for the second. In this
case we use coordinates $(t, z, w, \tau, \zeta, \eta)$ where we
write points $q' \in \Tscstar(X^\circ \times X^\circ \times \halfline)$
\begin{equation}
q' = \tau d(\frac1{t}) + \sum_{i=1}^n \zeta_i \frac{dz_i}{t} + \sum_{j=1}^n \eta_j \frac{dw_j}{t}.
\label{sc-coords}\end{equation}

Returning to the general situation, there is a contact structure defined at the boundary $\Tscstar[\pa X] X$ of $\Tscstar X$. It is defined by the contact one-form
\begin{equation}
\chi = \omega(x^2 \pa_x, \cdot) \restrictedto \{ x = 0 \},
\label{contactform}\end{equation}
where $\omega$ is the symplectic form on $T^* X$ (which is canonically isomorphic to $\Tscstar X$ over $X^\circ$). In local coordinates $\chi = \sum \mu_i dy_i  - d\nu$, so $\chi$ is clearly nondegenerate. A change of boundary defining function $x \to ax$ changes $\chi$ according to $\chi \to a \chi$, so the contact structure defined by $\chi$ is completely natural. A Legendrian submanifold of $\Tscstar[\pa X] X$ is defined, as usual,  to be a smooth submanifold of maximal dimension, namely $\dim X - 1$, such that the contact form $\chi$ vanishes on it. Any Legendrian submanifold $L$ has a local nondegenerate parametrization in a neighbourhood of any $q \in L$. This, by definition, is a function $\psi(y,v)$, with $v \in \RR^k$ for some $k \geq 0$, such that 
\begin{equation}\begin{gathered}
\text{ the differentials } d\big( \frac{\pa \psi}{\pa v_i} \big), \ i = 1, \dots, k \text{ are linearly independent whenever } d_v \psi = 0, 
 \\
\text{ and } L = \{ (y, d(\frac{\psi}{x})) \mid d_v \psi = 0 \} \text{ locally near } q.
\end{gathered}\label{ident}\end{equation}
The simplest situation is when $L$ is projectable in the sense that the projection $(y, \nu, \mu) \to y$ from $\Tscstar[\pa X] X$ to $\pa X$ restricts to a diffeomorphism from $L$ to $\pa X$. Then $y$ is a coordinate on $L$, so $\nu$ is given by a function $\psi(y)$ on $L$. In this case, no extra variables $v$ are required, and $L$ is given by $\{ (y, d(\psi/x)) \}$ locally.
If the kernel of the differential of the projection $(y, \nu, \mu) \to y$ restricted to $L$ has dimension $k$ then at least $k$ extra variables are required to locally parametrize $L$ near $q$.

A (half-density) Legendrian distribution of order $m$ associated to to $L$ is a half-density $u \alpha$, where $\alpha$ is a scattering half-density and $u$ is a finite sum of terms $\sum_i u_i + u_0$, where $u_0 \in \CdotI(X)$ and $u_i$ is given by an certain type of oscillatory integral associated to a local parametrization of $L$:
\begin{equation}
u_i=(2\pi)^{-k} x^{m-k/2+d/4} \int_K e^{i\psi(y,v)/x} a(y, v, x) \, dv.
\label{osc-int}\end{equation}
Here $a$ is a smooth function of $x, y, v$ with compact support and $\psi$ is a nondegenerate parametrization of $L$ on the support of $a$.  The set of
half-density Legendrian distributions of order $m$ associated with $L$ is denoted $I^m(X, L; \Omegasc^{1/2}(X))$, or just $I^m(L)$ when the space $X$ is understood (in this paper, Legendrian distributions will always be half-densities). 

Legendrian distributions have a well-defined symbol map $\sigma^m$
taking values in smooth sections of a line bundle $S^{[m]}(L)$ over
$L$. This bundle is given by $S^{[m]}(L) = \Omegasc^{1/2}(L) \otimes
\abs{N^* (\pa X)}^{m-d/4} \otimes M \otimes E$, where
$\Omega^{1/2}(L)$ denotes the half-density bundle over $L$, $N^*X$ is
the conormal bundle, $M$ is the Maslov bundle and $E$ is the bundle
described in \cite{MR2002i:58037}. To define the symbol, we choose $d-1$
functions $\lambda_j$ in $(y,v)$-space which together with $d_{v_i}
\psi$ give local coordinates in $(y, v)$-space. Then $\{ \lambda_j \}$ are local coordinates on $L$ via the identification \eqref{ident}. The
symbol is given, using the identification \eqref{ident} and up to Maslov factors, by
\begin{equation}
\sigma^m(u) = a(0, y, v) \abs{ \det \frac{\pa (\lambda, \pa_v \psi)}{\pa (y,v)} }^{-1/2} \abs{d\lambda}^{1/2} \restriction \{ \pa \psi/\pa v = 0 \}.
\label{legendriansymbol}
\end{equation}
The symbol of $u \in I^m(L)$ determines $u$ modulo $I^{m+1}(L)$. 
Consideration of how the symbol changes under changes of parametrization $\psi$, changes of coordinates $(x,y)$, and changes of coordinates $\lambda$, show that the symbol lives in the bundle $S^{[m]}(L)$ above \cite{MR2002i:58037}. If $L$ is locally projectable, then the situation simplifies. We may take coordinates $\lambda$ on $L$ to
be $y$, and the $v$ variables are absent, the determinant factor
above is $1$ and the symbol becomes
$$
a(0,y) \abs{dy}^{1/2}.
$$ 
The important bundle in the factorization of $S^{[m]}(L)$ here is  $\abs{N^* (\pa X)}^{m-d/4}$; in particular
this tells us that if we change boundary defining function from $x$ to
$x a^{-1}(y)$ then a symbol of order $m$ changes by a factor
$a^{m-d/4}$.  

We now recall the symbol calculus for a scattering pseudodifferential
operator $P \in \Psisc^{*,0}(X)$ acting on a Legendrian distribution
$u \in I^m(L)$. In fact, we only need to consider the case when the
symbol $p$ of $P$ vanishes identically on $L$, and in view of our
application in \S\ref{section:parametrix}, we use coordinates \eqref{sc-coords}
appropriate to the manifold $M = X^\circ \times X^\circ \times
\halfline$. Then the Hamilton vector field of the function $p$ (extended
into the interior of $\Tscstar M$ arbitrarily), $H_p$, vanishes to
first order at the boundary of $\Tscstar M$, so we define the rescaled
Hamilton vector field $\scH_p$ to be $x^{-1} H_p$ restricted to
$\Tscstar[\pa M] M$.  Then $Pu \in I^{m+1}(L)$ with symbol
\begin{equation}
\sigma^{m+1}(Pu) = \Big( -i {\mathcal{L}}_{\scH_p} - i (\half + m - \frac{d}{4})\frac{\pa p}{\pa
\tau} + p_{\sub} \Big) \sigma^m(u) \otimes \abs{dt}.
\label{transport}\end{equation}
Here $d=2n+1$ is the dimension of $M$, and
$\mathcal{L}_{H_p}$ denotes the Lie derivative acting on half-densities. The operation of
tensoring with $\abs{dt}$ is simply the natural isomorphism between
$S^{[m]}(L)$ and $S^{[m+1]}(L)$.

\subsection{Quadratic scattering structure}\label{qsc-str} We only touch on this very
briefly since we only make fleeting use of it in this paper. The
quadratic-scattering Lie algebra is defined by $\mathcal{V}_{\QSC} = x
\mathcal{V}_{\SC}(X)$. It is the space of smooth sections of a vector
bundle $\Tqsc X$ which is locally spanned (near the boundary) by $x^3
\pa_x$ and $x^2 \pa_{y_i}$. The dual bundle is denoted $\Tqscstar X$;
sections of this bundle are locally spanned over $\CI(X)$ by $dx/x^3$
and $dy_i/x^2$.  The vector fields in $\mathcal{V}_{\QSC} (X)$ lie in
a calculus of quadratic scattering pseudodifferential operators much
like the scattering calculus, and with a boundary symbol and
associated wavefront set lying in $\Tqscstar[\pa X] X.$  It is in terms
of this wavefront set that the propagation results of \cite{Wunsch1}
are couched, but it will not directly concern us here.

\subsection{Manifolds with corners with fibred boundaries}\label{subsec:fib-sc}
Above, we have looked at $X^\circ \times X^\circ \times \halfline$. However, in our analysis of the Schr\"odinger propagator $U(z, w,t)$, we wish to let the first variable $z$ approach the boundary of $X$. Hence we need to study the 
manifold
\begin{equation}
\XXt = X \times X^\circ \times \halfline.
\label{XXt}
\end{equation}
We denote the two boundary hypersurfaces of $\XXt$ at $t=0$ and at
$x=0$ by mf and sf, the `main face' and the `side face',
respectively, and we denote the corner $\mf \cap \tf$ by $K$. Thus
$\mf = X \times X^\circ$, $\tf = \pa X \times X^\circ \times \halfline$
and $K = \pa X \times X^\circ$.  We shall see that the geometry of this manifold
arising from the structure of the Schr\"odinger operator $D_t + H$ is
that of a manifold with corners with fibred boundaries; see \cite{MR2002i:58037}
and \cite{MR2001d:58034} for a discussion of the general situation. Here, we
restrict ourselves to describing the particular situation of $\XXt$,
which is considerably simpler.

We first describe the relevant fibration. Define $\phi : \tf \to \pa X \times X^\circ$ by projection off the last factor:
\begin{equation}
\phi(y, w, t) = (y, w) \in \pa X \times X^\circ.
\label{phi-def}\end{equation}
We then define a Lie algebra $\Vphi$ of smooth vector fields on $\XXt$ as follows: if $V$ is a smooth vector field on $\XXt$, then 
\begin{multline}
V \in \Vphi \text{ if and only if } V(xt) = O(x^2t^2), \ 
V \text{ is tangent to } \\ \text{ the fibration $\phi$ at $x=0$ and $V$ vanishes at $t=0$}. 
\end{multline}
This Lie algebra is independent of the choice of boundary defining function $x$. Corresponding to this Lie algebra is a vector bundle $\TphiXXt$ such that $\Vphi$ is canonically identified with the space of smooth sections of $\TphiXXt$. In local coordinates, this is spanned (near the boundary $x=0$) by 
\begin{equation}
xt \, \pa_{y_i}, \quad xt \, \pa_{w_i}, \quad t (t \pa_t - x \pa_x), \quad t x^2 \pa_x.
\label{sf-vfs}\end{equation}
Note that $t^2 \pa_t$ and $t x \pa_x$ are not separately in $\Vphi$, only their difference is. The dual bundle $\TstarphiXXt$ is spanned over $\CI(\XXt)$ by the one-forms
$$
d\big( \frac1{xt} \big) = -\frac{x dt + t dx}{x^2 t^2}, \quad d\big( \frac1{t} \big), \quad \frac{dy_i}{xt}, \quad \frac{dw_i}{xt}.
$$
Thus every point in $\TstarphiXXt$ near $\{ x = 0 \}$ may be expressed
$$
\sigma d(\frac{1}{xt}) + \tau d\big( \frac1{t} \big) + \overline{\mu} \frac{dy_i}{xt} + \overline{\xi} \frac{dw_i}{xt}
$$
($\overline{\mu}$ and $\overline{\xi}$ are used since we are reserving the symbols $\mu$ and $\xi$ for later use), which yields local coordinates $(x,y,w, t, \sigma, \overline{\mu}, \overline{\xi}, \tau)$ on $\TstarphiXXt$. 
A fibred-scattering half-density is defined to be a section of $\bigwedge^{2n+1}(\TstarphiXXt)$, $n = \dim X$, which is smooth and nonvanishing. Near the boundary, such a half-density has the form 
$$
a \Big| \frac{dx \, dy \, dw \, dt}{t^{2n+2}x^{2n+1}} \Big|^{1/2}, \text{ where } a \in \CI(\XXt) \text{ is everywhere nonzero.}
$$

On ${}^\phi T^* _{\mf} \XXt$, (that is, at $t=0$), $\omega(xt^2 \pa_t, \cdot)$ is a contact form (i.e.\ nondegenerate) for $x>0$, but it degenerates at $x=0$. This is evident from its local coordinate expression
$$
\omega(xt^2 \pa_t, \cdot) \restriction \{ t = 0 \} = d\sigma  + x d\tau + \overline{\mu} dy + \overline{\xi} dw.
$$
As described in \cite{MR2001d:58034}, although this is degenerate at ${}^\phi T^* _{K} \XXt$, that is, at ${}^\phi T^* _{\mf} \XXt \cap \{ x = 0 \}$, there is a natural fibration $\beta$ from ${}^\phi T^*_{K} \XXt$ to a bundle $\gamma$ over $K$, with a nondegenerate contact form on $\gamma$ such that $\omega(xt^2 \pa_t, \cdot) $ restricted to $\{ t= 0, x = 0 \}$ is the lift of the contact form from $\gamma$ to ${}^\phi T^*_{K} \XXt$. Here the bundle $\gamma$ is naturally isomorphic to $\Tscstar[\pa X\times X^\circ] (X\times X^\circ)$, the contact form is the natural contact form \eqref{contactform} on 
$\Tscstar[\pa X\times X^\circ] (X\times X^\circ)$ and the fibration $\beta$  is given by 
\begin{equation}
\beta(0,y,w,0, \sigma, \overline{\mu}, \overline{\xi}, \tau) = (y,w,\sigma, \overline{\mu}, \overline{\xi})
\label{beta-fib}\end{equation}
in the coordinates above. 

A fibred-scattering Legendrian submanifold $J$ of $\TstarphiXXt$, as defined in \cite{MR2001d:58034}, is a submanifold of ${}^\phi T^*_{\mf} \XXt$ which is Legendrian with respect to the contact form $\omega(xt^2 \pa_t, \cdot)$ for $x > 0$, which meets the boundary of ${}^\phi T^*_{\mf} \XXt$ at $x=0$ transversally, and such that the fibration $\beta$, restricted to $\pa J$, is a diffeomorphism from $\pa J$ to a Legendrian submanifold $G$ of $\Tscstar[\pa X\times X^\circ] (X\times X^\circ)$ (with respect to the natural contact structure on $\Tscstar[\pa X\times X^\circ] (X\times X^\circ)$). 

A local nondegenerate parametrization of $J$ near $q \in \pa J$, is a function $\psi(x,y,w,v)$,  with $v \in \RR^k$ for some $k \geq 0$, such that, locally near $q$,  
\begin{equation}\begin{gathered}
\text{ the differentials } d\big( \frac{\pa \psi}{\pa v_i} \big), \ 1 \leq i \leq k \text{ are linearly independent}
 \\ {\hskip 80pt}
\text{ and } J = \{ (y, d(\frac{\psi}{xt})) \mid d_v \psi = 0 \}. \end{gathered} \label{ident-fib}\end{equation}
It then follows that the function $\psi(0, y, w, v)$ parametrizes the Legendrian $G$. 

A Legendrian distributional half-density of order $(m,r)$ associated to $J$ is a half-density $u \tilde \alpha$, where $\tilde \alpha$ is a fibred-scattering half-density and $u$ is a finite sum of terms $\sum_j u_j + u_0 + u_f$, where $u_0$ is a Legendrian distribution of order $m$ associated to the interior of $J$ (where the structure is locally of scattering type, so this has already been defined), $u_f$ is a fibred Legendrian distribution of order $r$ associated to $G$, which is supported away from mf (defined below), and $u_j$ is given by an certain type of oscillatory integral associated to a local parametrization of $J$. Namely, $u_j$ is given by
\begin{equation}
(2\pi)^{-k} t^{m-k_j/2+(2n+1)/4} x^{r-k/2+(2n-1)/4} \int_E e^{i\psi_j(x,y,w,v)/xt} a_j(x,y, w,v, t) \, dv.
\label{sf-Leg}\end{equation}
Here $v \in \RR^{k_j}$, $E \subset \RR^{k_j}$ is compact, $a_j$ is a smooth function of $x, y, w,v$ with compact support and $\psi$ is a nondegenerate parametrization of $L$ on the support of $a$. 
Similarly, $u_f$ is given by a sum of terms of the form 
\begin{equation}\label{fibred-term}
(2\pi)^{-k} x^{r-k/2+(2n-1)/4} \int_E e^{i\psi_j(0,y,w,v)/xt} b_j(x,y, w,v, t) \, dv.
\end{equation}
with $b$ smooth and supported away from $t=0$. It is straightforward to check that \eqref{sf-Leg} and \eqref{fibred-term} are compatible. Note that, for fixed $t > 0$, the expression \eqref{fibred-term} is a Legendrian distribution of order $r-1/4$ associated to the Legendre submanifold $t^{-1} G$, where $t^{-1}$ acts by scaling in the fibre variables. 

Fibred-scattering Legendrian distributions of order $(m,r)$ on $L$ have a well-defined symbol map $\sigma^m$ taking values in $x^{r-m} \CI(L;S^{[m]}(L))$ over $L$, which is defined by continuity from the symbol map in the interior of $L$. (This makes sense because the fibred-scattering structure at mf and away from sf is just the scattering structure.) The bundle $S^{[m]}(L)$ in the fibred-scattering setting is given by
$S^{[m]}(L) = \Omega_b^{1/2}(L)
\otimes \abs{N^* \mf}^{m-N/4} \otimes \abs{N^* \tf}^{m-N/4} \otimes M \otimes E$, where $\Omega_b^{1/2}(L)$ denotes the $b$-half-density bundle over $L$; a smooth nonzero section near the boundary is $|d\lambda dx/x|^{1/2}$, where $\lambda$ are coordinates on $\pa L$, extended into the interior of $L$. The symbol determines $u \in I^{m,r}(L)$ up to an element of $I^{m+1,r}(L)$.


\section{Geometry of the time-dependent Schr\"odinger operator}\label{section:flowout}
\subsection{Flowout from the diagonal}
Let $\so=t^2 (D_t + H)$. We consider this operator as an element of order $(2,0)$ of the scattering calculus for the manifold with boundary $X^\circ \times X^\circ \times \halfline $. We shall use coordinates $z, w, t, \zeta, \eta, \tau$ in the scattering cotangent bundle, as in \eqref{sc-coords}.  The boundary symbol of  $\so$ at $t=0$ is
$$
p(z, w, \zeta, \eta, \tau) = -\tau + \frac1{2} \abs{\zeta}_g^2.
$$
The rescaled Hamiltonian flow is given by 
\begin{equation}\begin{gathered}\begin{aligned}
\dot z &= \h \frac{\partial \abs{\zeta}^2_g}{\partial \zeta} \\ 
\dot w &= 0 \\
\dot t &=t\\ 
\end{aligned}\end{gathered}\qquad \qquad \begin{gathered}\begin{aligned}
\dot \zeta &= \zeta - \h\frac{\pa \abs{\zeta}^2_g}{\pa z}\\ 
\dot \eta &= \eta \\
\dot \tau &= \abs{\zeta}^2_g .
\end{aligned}\end{gathered}\label{Hflow-int}\end{equation}
The free Schr\"odinger propagator $c_n t^{-n/2} e^{i|z-w|^2/2t}$ is a Legendrian distribution parametrized by the phase function $|z-w|^2/2$. We look for an analogue in the general case. 

\begin{lemma}\label{Leg}
There exists a smooth Legendre submanifold $L$ of $\Tscstar[ \{ t = 0 \} ] X^\circ \times X^\circ \times \halfline$ contained in $\{ p = 0 \}$  that is parametrized by the phase function 
$\Phi(z,w) = d(z,w)^2/2$ near the diagonal $z=w$. 
\end{lemma}

\begin{proof} Since $p(z, w, d_z \Phi, d_w \Phi, \Phi) = 0$, the function $\Phi$ 
parametrizes a Legendrian contained in $\{ p = 0 \}$ in the region where $\Phi$ is smooth, that is, within the injectivity radius. We define $L$ to be the flowout along bicharacteristics from this region. $L$ is given as a set by 
\begin{equation}
(z, \hat \zeta) = \exp_{sg/2} (w, \hat \eta_0), \ \abs{\eta} = \abs{\zeta} = s, \ \tau = \frac{ s^2}{2}, \ \hat \eta = -\hat \eta_0, \ s \in [0, \infty), \ (w, \hat \eta_0) \in S^* X^\circ .
\label{geod-flow}\end{equation}
Since the Hamilton vector field never vanishes for $s > 0$, and since $\tau \to \infty$ as $s \to \infty$, the entire flowout is smooth. 
\end{proof}

\subsection{Behaviour of $L$ near the boundary}
We next analyze the behaviour of $L$ near the boundary of $X$ in the first factor. To do this, we introduce the Lie algebra of vector fields given by the $C^\infty(X \times X \times \halfline)$-span of the vector fields
$$
t^2 x^2 \pa_t, \ t x^3 \pa_x, \ t x^2 \pa_y,\  t x^2 \pa_w.
$$ 
There is a vector bundle of
which these are smooth sections forming a basis at every
point. Let $B$ denote the dual of this bundle; then a basis of the smooth sections of $B$ is given locally near the boundary of $X$ by
$$
\frac{dt}{t^2 x^2}, \frac{dx}{x^3 t}, \frac{dy_i}{x^2 t}, \frac{dw_j}{x^2 t}.
$$
Hence an arbitrary $q \in B$ may be written\footnote{The variable $\mu$ here is not the same as the $\mu$ from Sections~\ref{section:intro} and \ref{section:definitions}.}
\begin{equation}
q = \frac{\kappa}{x^2} d( \frac1{t}) + \frac{\lambda \, dx}{x^3 t} + \frac{\mu \, dy}{x^2 t} +
\frac{\xi \, dw}{x^2t}.\label{coords}
\end{equation}
Then $x, y, w, t, \lambda, \mu, \xi, \kappa$ are smooth coordinates on
$B$, up to the corner $t = x = 0$.  We remark that the bundle $B$ is a
scaled cotangent bundle mixing scattering behaviour at the face $t=0$ with
quadratic scattering behaviour at $x=0.$

There is a natural vector bundle isomorphism from $\Tscstar (X^\circ\times X^\circ \times \halfline)$ to $B$ over $X^\circ$ given by
\begin{equation}
(z, w, t, \zeta, \eta, \tau) \to (z, w, t, x \zeta, x\eta, x^2 \tau)
\label{B-iso}\end{equation}
in terms of the coordinates given above. Thus, the map rescales the
cotangent variables by a power of $x$ as we approach $x=0$ --- linear
rescaling for the spatial cotangent variables and quadratic for the
temporal cotangent variable. This counteracts the growth of these
variables as the boundary $x=0$ is approached, since we see from
\eqref{geod-flow} that $\zeta, \eta$ grow linearly and $\tau$
quadratically in $s$ (which expected to be asymptotic to $x^{-1}$ as $s \to
\infty$). Thus, if we map $L$ into $B$ via the isomorphism
\eqref{B-iso}, we can expect it to remain `bounded' at $x=0$. By abuse
of notation, we denote the image of $L$ under \eqref{B-iso} also by
$L$. The Legendrian $L$, now regarded as a submanifold of $B$, is
generated by Hamiltonian flow of $x^2 \so$ starting from the diagonal.

Owing to the special structure \eqref{better} of the scattering metric $g$, 
the boundary symbol of $x^{2} \so$ equals
\begin{equation}
-\kappa + \h (\lambda^2+h^{ij}\mu_i \mu_j) 
\label{corner-symbol}\end{equation}
where $h^{ij}$ is a smooth function of $x$ and $y$. This
yields (rescaled) Hamilton equations in $B$
\begin{equation}\begin{gathered}\begin{aligned}
\dot x &= \lambda x \\
\dot y_i &= h^{ij}\mu_j  \\
\dot w &= 0\\
\dot t &= t 
\end{aligned}\end{gathered} \qquad \qquad \begin{gathered}\begin{aligned}
\dot \lambda &= \lambda+ 2\kappa -2h^{ij}\mu_i \mu_j - \frac{x}{2} \frac{\pa h^{ij}}{\pa x} \mu_i \mu_j  \\
\dot \mu_i &= (1 + 2\lambda) \mu_i-\h \pa_{y_i} h^{jk}\mu_j \mu_k \\
\dot {\xi} &= (1+2\lambda)\xi \\
\dot \kappa &= \lambda^2+ 2\lambda \kappa + h^{ij}\mu_i \mu_j . \end{aligned}\end{gathered}\label{hamflow1}\end{equation}

The next lemma justifies our hope that the Legendrian $L$ stays bounded
 in $B$ as we approach the boundary, $x \to 0$. 

\begin{lemma}\label{asympt-ray}
Along every non-trapped bicharacteristic ray inside the Legendrian $L$, we have $x \to 0$, $\xi \to 0$, $\lambda \to -1$, $\mu \to 0$.
\label{lemma:attractor}
\end{lemma}

\begin{proof}
Along every nontrapped ray, we have $x \to 0$. Hence, there are points along the ray where $x$ is arbitrarily small and where $\dot x < 0$. Since $\dot x = \lambda x$, there are points along the ray where $x$ is arbitrarily small and where $\lambda < 0$. 

Now we look at the equation for $\lambda$. Since the ray lies inside the characteristic variety of $x^2 \so$, we can substitute $2\kappa = \lambda^2 + h^{ij} \mu_i \mu_j$ on the right hand side of \eqref{hamflow1}. Hence
\begin{equation}
\dot \lambda = \lambda + \lambda^2 - h^{ij} \mu_i \mu_j - \frac{x}{2} \frac{\pa h^{ij}}{\pa x} \mu_i \mu_j  .
\label{ldot}\end{equation}
Since $h^{ij}(x)$ is positive definite for small $x$, uniformly over $\pa X$ (since $\pa X$ is compact), 
$$
-h^{ij} \mu_i \mu_j - \frac{x}{2} \frac{\pa h^{ij}}{\pa x} \mu_i \mu_j \leq 0 \text{ for } x < \epsilon. 
$$
Thus, starting from a point on the ray where $x < \epsilon$ and $\lambda < 0$, we have
$$
\dot x = \lambda x, \quad \dot \lambda \leq \lambda + \lambda^2 \text{ if } x < \epsilon. 
$$ These equations imply that $\lambda$ remains negative and that $x$
is decreasing from this point forwards on the ray. Moreover, since
$\lambda + \lambda^2 < 0$ for $-1 < \lambda < 0$, we see that $\limsup
\lambda \leq -1$. In particular, from some point on, $\lambda <
-3/4$. This implies that $\xi \to 0$ along the ray.

Now consider the equation for $\mu$. It implies that 
$$
\big( h^{ij} \mu_i \mu_j \big)^{{\bf \cdot}} = 2 (1 + 2\lambda) h^{ij} \mu_i \mu_j + x \lambda \frac{\pa h^{ij}}{\pa x} \mu_i \mu_j.
$$ Since eventually $\lambda < -3/4$ and $x < \epsilon$, positivity of
$h^{ij}$ implies that the right hand side is $\leq -1/2 h^{ij} \mu_i
\mu_j$, say, from some point on, so $\mu \to 0$. Returning to the
equation \eqref{ldot} for $\lambda$, we see that $\mu \to 0$ implies
that $\lambda \to -1$, since \eqref{ldot} with the $\mu$ terms removed
has an attracting fixed point at $\lambda = -1$.
\end{proof}

\subsection{Smoothness at the boundary}
Throughout this subsection, we assume that the metric $g$ is nontrapping. 

Let $\Sigma(\so)$ be the zero set of the boundary symbol of $x^2 \so$ on $B \cap \{ t = 0 \}$.  We identify this set with $\Tqscstar(X\times
X^\circ)$ by mapping
$$
\frac{\kappa}{x^2} d( \frac1{t}) + \frac{\lambda \, dx}{x^3 t} + \frac{\mu \, dy}{x^2 t} +
\frac{\xi \, dw}{x^2t} \qquad \text{ to } \qquad  \frac{\lambda \, dx}{x^3} + \frac{\mu \, dy}{x^2} +
\frac{\xi \, dw}{x^2}.
$$
In terms of coordinates, we map $(x, y, w, \kappa, \lambda, \mu, \xi)$ to $(x, y, w, \lambda, \mu, \xi)$. This is an isomorphism on $\Sigma(\so)$ since $\kappa$ is given by $\abs{(\lambda, \mu)}_{g_{x,y}}^2 /2$ on $\Sigma(\so)$. 
The image of $L$ under this identification is then a {\it Lagrangian} in $\Tqscstar(X\times X^\circ)$ (which we shall denote $\Lt$). We emphasize that this Lagrangian lies over the interior of $X \times X^\circ$, not at the boundary. In terms of parametrizations, if $\Phi/t$ is a local parametrization for $L$, then $\Phi$ is a local parametrization of $\Lt$. 

It turns out that $\Lt$  is not smooth at $x=0$. To recover
smoothness, we shall blow up the submanifold 
\begin{equation}
S = \{ x=0, \lambda = -1, \mu = 0,
\xi=0 \}
\label{S}\end{equation}
where $\Lt$ meets the boundary (by Lemma~\ref{asympt-ray}). 
We shall use Melrose's notation and terminology; thus, we denote the manifold with corners obtained by blowing up the submanifold $S$ by 
$$
\big[ \Tqscstar(X\times X^\circ); S \big];
$$ explicitly, this is the manifold with corners obtained by removing $S$
and replacing it with its inward-pointing spherical normal bundle, with a
natural $\CI$ structure (see \cite{M}). This inward-pointing spherical
normal bundle becomes a new boundary hypersurface of the blown-up space
which we refer to as the `front face' of the blowup.  The lift to this
blown-up space of the original boundary hypersurface $\pa X \times X^\circ$
will be denoted $\widetilde \tf$.

\begin{lemma}\label{Lt}
\label{smooth} Let $(X,g)$ be a manifold with boundary, with a nontrapping scattering metric. Let $S \subset \Tqscstar(X
\times X^\circ)$ be given by \eqref{S}, where $x$ is a boundary defining function  such that $|dx/x^2|_{g} = 1$ at $\pa X$.  Then
\begin{equation}
\XtS \equiv \big[ \Tqscstar(X\times X^\circ); S \big] \setminus \widetilde \tf
\label{Xts}
\end{equation}
is an affine bundle over $X \times X^\circ$, isomorphic (but not
naturally isomorphic) to the scattering cotangent bundle $\Tscstar(X \times
X^\circ)$.  The Lagrangian $\Lt$ lifts to the space \eqref{Xts} to be a smooth manifold with boundary, such that $\pa \Lt$ is contained in the boundary of
$\Tscstar(X \times X^\circ)$ and such that $\Lt$ is transversal to the
boundary there.
\end{lemma}

\begin{proof}
Performing the blowup amounts to introducing the coordinates  on $\XtS$ in a neighbourhood of the front face. The induced map from $\XtS$ to $X \times X^\circ$ is a
fibration in which each fibre is diffeomorphic to $\RR^{2n}$. To specify the affine structure we introduce the coordinates 
$\Lambda = (\lambda+1)/x$, $M = \mu/x$, $\Xi = \xi/x$ which are smooth on $\XtS$ up to $x=0$. The affine
structure is defined by deeming $\Lambda, M$ and $\Xi$ to be affine
coordinates on each fibre. This is well defined since a smooth change of
variables in $(x, y, w, \lambda, \mu, \xi)$ which is linear on each fibre
of $\Tqscstar(X\times X^\circ)$ induces a smooth change of variables in
$(x, y, w, \Lambda, M, \Xi)$ which is {\it affine} on each fibre. (Observe
that a coordinate change such as $x \to x + \alpha x^2$ induces a
translation $\Lambda \to \Lambda + \alpha + O(x)$, so the coordinate
changes are not in general linear on the fibres.) In terms of the new
coordinates, a point
$$
q = \frac{\lambda \, dx}{x^3} + \frac{\mu \cdot dy}{x^2} + \frac{\xi \cdot dw}{x^2} \in \Tqscstar(X \times X^\circ)
$$
becomes
\begin{equation}
q = -\frac{dx}{x^3}+ \Lambda\frac{dx}{x^2}+ M \cdot \frac {dy}{x} + \Xi \cdot \frac{dw}{x} .
\label{sc-bundle}\end{equation}
Mapping $q$ to $\Lambda dx/x^2 + M dy/x + \Xi dx/x$ defines an affine
bundle isomorphism between the blown-up space $\XtS$ defined in \eqref{Xts}
and $\Tscstar(X\times X^\circ)$.  It is not, however, a canonical
isomorphism --- this is clear, since the scattering cotangent bundle has a
linear structure which $\XtS$ lacks. 
We remark that if we performed the blowup at the zero section of the
quadratic-scattering cotangent bundle, that is, at $\lambda = 0$ instead of
$\lambda = -1$, then the resulting space as in \eqref{Xts} {\it would} be naturally isomorphic to the
scattering cotangent bundle.

Next, we show that $\Lt$ has the regularity specified in the statement of the lemma. 
Our method of proof is to lift the vector field $V$ in \eqref{hamflow1} to the space $\big[ \Tqscstar(X\times X^\circ); S \big]$, and divide by a boundary defining function for the front face, obtaining a new vector field $W$. Clearly integral curves of $V$ are the same as integral curves of $W$. We claim that $W$ is a smooth vector field on $\big[ \Tqscstar(X\times X^\circ); S \big]$ which is tangent to $\widetilde \tf$ but transverse to the front face. Consequently, every integral curve of $W$ must meet the boundary of $\big[ \Tqscstar(X\times X^\circ); S \big]$ in the interior of the front face. Moreover, the fact that $W$ is smooth and nonvanishing at the front face implies that the flowout $\Lt$ is smooth up to the front face and meets it transversally. Hence to prove the lemma, it is sufficient to prove the claim above. 

Proving the claim is simply a matter of computing the vector field \eqref{hamflow1} in various coordinate systems valid in different regions of the blown-up space $\big[ \Tqscstar(X\times X^\circ); S \big]$. In a neighbourhood of any point in the interior of the front face, we may use coordinates $\Lambda, M, \Xi$ defined above. In terms of these coordinates, the vector field takes the form 
\begin{multline}
 V=(-x+\Lambda x^2) \pa_x + x h^{ij} M_j \pa_{y_i} +x (- h^{ij} M_i
 M_j + \frac{x}{2} \frac{\pa h^{ij}}{\pa x} M_i M_j ) \pa_\Lambda \\ +
 x (\Lambda M_i- \frac{1}{2} \frac{\pa h^{jk}}{\pa y_i } M_j M_k)
 \pa_{M_i} + x \Lambda \Xi \cdot \pa_{\Xi}.
\label{hamflow:blownup}
\end{multline}
Dividing by $x$, which is locally a boundary defining function for the front face in this region, shows that $W = V/x$ is a smooth vector field in this region which is transverse to the boundary. 

Near the corner of $\big[ \Tqscstar(X\times X^\circ); S \big]$, i.e.\
near the boundary of the front face, we need several different
coordinate patches. Let $\rho = \sqrt{x^2 + (\lambda + 1)^2 + |\xi|^2
  + |\mu|^2}$, and let $s = \lambda+1$. When $s \geq \rho/4$, then
coordinates $s_1, r_1 = x/s_1, M_1 = \mu/s_1, \Xi_1 = \xi/s_1$ may be
used. When $|\mu| \geq \rho/4$, then valid coordinates are $s_2 =
|\mu|, r_2 = x/s_2, \Lambda_2 = (\lambda+1)/s_2, \hat \mu = \mu/s_2,
\Xi_2 = \xi/s_2$ may be used, and when $|\xi| \geq \rho/4$, then valid
coordinates are $s_3 = |\xi|, r_3 = x/s_3, \Lambda_3 =
(\lambda+1)/s_3, M_3 = \mu/s_3$, $\hat \xi = \xi/s_3$. The union of
these coordinate patches covers a neighbourhood of the corner of
$\big[ \Tqscstar(X\times X^\circ); S \big]$.

In the first set of coordinates, $(s_1, r_1, M_1, \Xi_1)$, the vector field $V$ takes the form (where the subscript $1$ is dropped for notational convenience)
\begin{multline}
V = (-s + s^2 k) \pa_s - rsk \pa_r + sh^{ij} M_j \pa_{y_i} \\  + s((1 - sk) M_i - \frac{1}{2} \frac{\pa h^{jk}}{\pa y_i} M_j M_k ) \pa_{M_i} + s(1 - sk) \Xi_i \pa_{\Xi_i}
\end{multline}
The vector field $W = V/s$ is tangent to $\widetilde \tf$ and
transverse to the front face in the region of validity of these
coordinates. The calculation for the second and third set of
coordinates is similar, and is left to the reader. 
\end{proof}

Now let $G$ be the boundary of $\Lt$: $G=\pa \Lt\subset \Tscstar[\pa X \times X^\circ] (X\times X^\circ)$.
\begin{lemma}\label{G} The submanifold $G$ of $\Tscstar[\pa X \times
    X^\circ] (X\times X^\circ)$ is Legendrian. 
\end{lemma}

\begin{proof} 
 The contact form on $\Tscstar[\pa X \times X^\circ] (X\times X^\circ)$ is
given by $\omega(x^2 \pa_x, \cdot),$ where the symplectic form
$\omega$ is now given by
\begin{equation}
\omega=\frac{d\Lambda\wedge dx}{x^2}+ \frac{dM\wedge dy}{x}+ \frac{d\Xi \wedge dw}{x}- \frac{M dx \wedge dy}{x^2} - \frac{ \Xi dx \wedge dw}{x^2}
\label{omega}\end{equation}
which we obtain by taking the differential of \eqref{sc-bundle}.
Since $V$ in \eqref{hamflow:blownup} is of the form $-x \pa_x + x
V''$, where $V''$ is tangent to the boundary of $\Tscstar (X\times
X^\circ)$, and $x\omega$ by \eqref{omega} is nondegenerate on vector
fields tangent to the boundary, uniformly up to the boundary, we find
that the contraction of $\omega$ and $x^2 \pa_x$ satisfies
$$
\omega(x^2 \pa_x, \cdot) =  -x\omega(V,\cdot) + x\alpha
$$
where $\alpha$ is a smooth one-form on $\Tscstar
(X\times X^\circ)$. Let $W$ be a vector tangent to $G$. Then 
$$
\omega(x^2 \pa_x, W) =  -x\omega(V,W) + x\alpha(W),
$$
The first term vanishes since $G \subset \tilde L$, both $V$ and $W$
are tangent to $\tilde L$, and $\tilde L$ is Lagrangian. The second term vanishes at 
$x=0$. This proves that $G$ is Legendrian.
\end{proof}

\subsection{Asymptotics of geodesic flow and the Legendrian $G$}
The variables $\Lambda$ and $M$ have geometric meaning on the Legendrian $G$ in terms of the asymptotic behaviour of geodesics. Notice that, by \eqref{geod-flow}, the value of $\tau$ on the Legendrian $L$ is given by $d(z, w)^2/2$. 
Thus, $\kappa$ is given by 
$$
\kappa = x^2\tau = x^2 \, \frac{d(z,w)^2}{2}.
$$
Near the boundary $G$ of the Legendrian, we have by \eqref{corner-symbol}
$$
\kappa = \frac1{2} (\lambda^2 + h^{ij}\mu_i \mu_j) = \frac1{2} \big( (-1 + x\Lambda)^2 + x^2 h^{ij} M_i M_j \big).
$$
This implies that
\begin{equation}
d(z,w) = \frac1{x} - \Lambda + O(x).
\label{Lambda-soj}\end{equation} 
Thus we interpret $-\Lambda$ on $G$ as a `sojourn time', the time that the
geodesic spends in the interior of the manifold before emerging into
the conic region. 

\begin{remark}
It is instructive to note that under a change of boundary defining function
$x \to \tilde x = x+ \alpha x^2$ which preserves the normal form
\eqref{better} we find
that the sojourn time changes as follows: letting $t$ be a unit speed
parameter along the geodesic flow, we compute
$$
\lim t-\frac 1{\tilde{x}} =\lim t- \frac 1x \frac x{\tilde x} = -\alpha + \lim
t-\frac 1x;
$$ hence we easily verify that this `sojourn time,' while ill-defined as a
function, does transform as a coordinate on the affine bundle
$\XtS\restrictedto_{x=0},$ as we know it must.
\end{remark}

To interpret the coordinate $M$, we notice that according to \eqref{hamflow:blownup} we have
\begin{equation} \begin{gathered} \begin{aligned}
\dot y_i &= h^{ij} \mu_j = x h^{ij} M_j , \\
\dot x &= -x + O(x^2),
\end{aligned}
\end{gathered} \quad \text{ which implies } dy_i/dx = -h^{ij} M_j+O(x).
\label{MMM}\end{equation}
Hence $M$ measures the
\emph{angle} at which the geodesic strikes the boundary, or more
geometrically, $M$ is a coordinate on the pencil of `asymptotically
parallel' rays with a fixed final direction $y$.  The coordinate $y$,
of course, is just the limiting direction of the geodesic, while $w$
and $-\hat \Xi$, which are constant under the flow, give the initial
conditions of the geodesic.

\subsection{Fibred scattering structure}

There is yet another way of viewing the Legendrian $L$, in terms of the fibred-scattering structure  of $\XXt$ described in Section~\ref{subsec:fib-sc}. We show that one can map the Legendrian $L \subset B$ to a fibred-scattering Legendrian $\Lphi$ of $\TstarphiXXt$. 
 To define this map, we use local coordinates $(x, y, w, t)$. Let \begin{equation}
q = \frac{\kappa}{x^2} d( \frac1{t}) + \frac{\lambda \, dx}{x^3 t} + \frac{\mu \, dy}{x^2 t} +
\frac{\xi \, dw}{x^2t},
\label{q}\end{equation}
be a point in $B$, using coordinates as in \eqref{coords}. We map $q\in L \subset B$ to 
\begin{equation}
q \mapsto q + \frac1{2} \frac{dt}{x^2 t^2} + \frac{dx}{x^3t} = q - \frac1{2} d(\frac1{x^2t}) \in \TstarphiXXt.
\label{add}\end{equation}
\begin{lemma}
Assume that the metric $g$ is nontrapping. Then the image of \eqref{add} is a smooth fibred-scattering Legendrian submanifold
$\Lphi$ of $\TstarphiXXt$.
\label{lemma:fibredLeg}
\end{lemma}

{\it Remark. }
Notice that the terms $dt/t^2 x^2$ and $dx/x^3t$ in \eqref{q} are too `big' for the fibred scattering space $\TstarphiXXt$. The addition of terms as in \eqref{add} cancels the big part of these terms and yields something which remains `bounded' in $\TstarphiXXt$ up to $x=0$. 

\begin{proof}
We use the blowup coordinates $\Lambda, M, \Xi$ as before. Thus,
$$\begin{aligned}
\lambda &= -1 + x\Lambda \\
\mu &= xM \\
\xi &= x\Xi \\
\kappa &= \half \big( \lambda^2 + h^{ij}\mu_i \mu_j \big) = -\half - x \Lambda + x^2 (\Lambda^2 + h^{ij}M_i M_j ) 
\end{aligned}$$
Now we add $\half dt/x^2 t^2 + dx/x^3 t$ to $q$ in \eqref{q}, and get
$$
q \mapsto \Lambda \frac{x dt + t dx}{x^2 t^2} + (\Lambda^2 + h^{ij}M_i M_j ) 
\frac{dt}{t^2} + M \frac{dy}{xt} + \Xi \frac{dw}{xt}.
$$
Thus there is a diffeomorphism between $\tilde L$ and $\Lphi$, given in local coordinates by
$$
(x, y, w, \Lambda, M, \Xi) \mapsto (x, y, w, \Lambda, \Lambda^2 + h^{ij}M_i M_j, M, \Xi),
$$
so $\Lphi$ is smooth up to the boundary and meets it transversally. 
The submanifold $\Lphi$ is a Legendrian distribution. Indeed, if $\Psi$ locally parametrizes $L$ away from the boundary of $B$, then $\Psi - 1/(2x^2 t)$ locally parametrizes $\Lphi$. Finally, under our identifications, the image of the boundary of $\Lphi$ under $\beta$ is identical with $G$ from Lemma~\ref{G}, which is a Legendrian submanifold. This proves that $\Lphi$ is a fibred-scattering Legendre submanifold. 
\end{proof}

\subsection{Trapped rays} If the metric $g$ has trapped rays, then we need to localize the constructions of this section in the nontrapping region. We define $L$ to be as in Lemma~\ref{Leg}, but taking only those $(w, \hat \eta) \in \nft$, rather than all $(w, \hat \eta) \in S^* X^\circ$. We then define $\Lt, \Lphi$ and $G$ in terms of $L$ as before. Then Lemmas~\ref{Lt}, \ref{G} and \ref{lemma:fibredLeg} continue to hold except we lose smoothness of $L$, $\Lt$ and $\Lphi$ at the set corresponding to $s=0$ in \eqref{geod-flow}. Moreover, if we localize to a compact subset of $\nft$, then the corresponding pieces of $L$, $\Lt$ and $\Lphi$ are compact.

\begin{lemma} The map $S_f$ defined by
\begin{equation}
S_f(w, \hat \Xi) = (y, \Lambda, M) \Longleftrightarrow
  (y,w,\Lambda, M, -\Xi) \in G, \quad \frac{\Xi}{|\Xi|} = \hat \Xi
\label{Sdefn}\end{equation}
is a contact diffeomorphism
$
S^* X^\circ\supset\nft \mapsto \Tscstar[\pa X] (X)
$
satisfying \eqref{soj-defn} and \eqref{pencil-coord}. 
Similarly, $S_b(w, \hat \Xi) = -S_f(w, -\hat \Xi)$ is a contact diffeomorphism.
\label{lemma:contactomorphism}
\end{lemma}

\begin{proof}
For each $(w, \hat \Xi) \in \nft$ there is a unique point $(y, w, -\Lambda, -M, \Xi) \in G$ found by following the bicharacteristic emanating from $(w, \hat \Xi)$ until it hits the front face of $\XtS$; uniqueness and smoothness are a consequence of the transversality of the vector field $W$ in the proof of Lemma~\ref{Lt}. 
Invertibility of the map $S_f$ follows from the previous discussion. In fact, given $(y, \Lambda, M) \in \Tscstar[\pa X](X)$, $y$ determines a pencil of geodesics with asymptotic direction $y$, $M$ picks out a unique geodesic within the pencil and then $\Lambda$ indicates how far along the geodesic one must travel to get to the initial point $w$. The direction $\Xi$ at $w$ is the direction of that geodesic at $w$. Hence for each $(y, \Lambda, M)$ there is a unique $(w, \hat \Xi)$ with $S_f(w, \hat \Xi) = (y, \Lambda, M)$. 
Comparison between \eqref{numu} and \eqref{sc-bundle} shows that $\Lambda = -\nu$ when interpreted as a coordinate on $\Tscstar X$, hence \eqref{Lambda-soj} implies \eqref{soj-defn}, while \eqref{MMM} implies \eqref{pencil-coord}. 

The fact that $S_f$ is contact follows directly from the definition \eqref{Sdefn} and the fact that $G$ is Legendrian (see \cite{GS}, p149, for a discussion in the symplectic case). 
\end{proof}


\section{Propagator (nontrapping case)}\label{section:parametrix}
In this section, we shall construct a parametrix
$\cU(z, w,t)$ for the Schr\"odinger propagator, where we restrict  the second variable
$w$ to some arbitrary open set $\cG$ with compact closure in
$X^\circ$. We require that $\cU$ solves
\begin{equation}
\so \cU (z,w,t) \equiv t^2(D_t + H_z) \cU(z, w,t) = \cE(z, w,t)\in t^N x^N \CI(X \times \mathcal{G} \times \halfline) \text{ for all } N,
\label{eqn}\end{equation}
or in other words, the error term $\cE$ is smooth and rapidly decreasing both as $t \to 0$ and as the $z$ variable tends to infinity. The restriction to $t >
0$ is only for convenience. In addition $\cU$ should satisfy the initial
condition
\begin{equation}
\cU(\cdot, w, t) \to \delta_w \text{ in } C^{-\infty}(X) \text{ as } t \to
0, \ \text{ for all } w \in \cG.
\label{init-cond}\end{equation}
In this section, we deal with the case in which the metric is globally nontrapping; we sketch the changes necessary to localize the construction in
$\nft$ in the following section.

The parametrix will be a sum of four terms, $\cU = \cU_1 + \cU_2 + \cU_3 +
\cU_4$. Correspondingly, the construction is divided into four steps. We
first construct a $\cU_1$ supported near $t=0$ and near the diagonal,
satisfying the initial condition \eqref{init-cond} and solving
\begin{equation}
\so \cU_1 = \cE_1 + \cR_1,
\label{eqn1}\end{equation}
where $\cE_1$ satisfies the condition for the error term in \eqref{eqn},
and $\cR_1( z, w,t)$ is supported in $\iota(X,g)/4 \leq d(z, w) \leq
\iota(X,g)/2$, where $\iota$ is the injectivity radius of $(X,g)$.  In the second step, we solve
\begin{equation}
\so \cU_2 = -\cR_1 + \cE_2 + \cR_2,
\end{equation}
with zero initial conditions, where $\cE_2$ satisfies the condition for the error term in \eqref{eqn}, and with $\cR_2(z, w,t)$ supported in $d(z,w) \geq \iota(X,g)/4$, $t \leq 1$ and with $x \geq 2\epsilon$. (Here, and henceforth, we use $(x,y)$ as local coordinates for the $z$ variable when it is close to the boundary.) In the third step, we solve 
\begin{equation}
\so \cU_3 = -\cR_2 + \cE_3 + \cR_3,
\end{equation}
with zero initial conditions, where $\cE_3$ satisfies the condition for the error term in \eqref{eqn}, and with $\cR_3(z, w, t)$ supported in $x \leq 4\epsilon$. Finally, we solve 
\begin{equation}
\so \cU_4 = -\cR_3 + \cE_4,
\end{equation}
with zero initial conditions, where $\cE_4$ satisfies the condition for the error term in \eqref{eqn}. Clearly the sum $\cU_1 + \cU_2 + \cU_3 + \cU_4$ satisfies the conditions for a parametrix.

\

\subsection{Step 1 --- near the diagonal}\label{Step1} 
We start by constructing a formal solution to
$$
\so \tilde \cU_1(z, w, t) = t^2 (D_t + H) \tilde \cU_1(z, w, t) = 0, \quad t > 0,
$$
near $z = w$, with initial condition a delta function $\delta_w(z)$. 
Our ansatz is
\begin{equation}
\tilde \cU_1 = t^{-n/2} e^{i\Phi(z,w)/t} \sum_{j=0}^\infty t^j a_j(z, w).
\end{equation}
By a formal solution
we mean that each coefficient of $t^j$ vanishes (near $z=w$) after applying the operator.
Applying the operator, we find
\begin{equation}\begin{gathered}
t^{n/2}\so \tilde \cU_1 = e^{i\Phi(z,w)/t} (-\Phi + g(\nabla_z \Phi, \nabla_z \Phi)) \sum_{j=0}^\infty t^j a_j \\
+ \,  it \, \Big( -g(\nabla \Phi,\nabla) a_0 + \frac{n}{2} a_0 + \half (\Delta \Phi) a_0 \Big)
\\
+ \sum_{j=1}^\infty it^{j+1} \Big( -g(\nabla \Phi,\nabla) a_j + (\frac{n}{2} - j) a_j + 
\half (\Delta \Phi) a_j + \frac{i}{2} \Delta a_{j-1} - i V a_{j-1} \Big)
\end{gathered}\label{formal-1}\end{equation}  
This gives us a sequence of equations to be solved, one for each
power of $t$. The first equation is the {\it eikonal equation}, $-\Phi +
g(\nabla_z \Phi, \nabla_z \Phi) = 0$.  This has an exact solution
\begin{equation}
\Phi = \half d(z,w)^2, 
\label{eik}\end{equation}
which is smooth when $d(z,w)$ is smaller than the injectivity radius of
$(X,g)$.  Motivated by the form of the free propagator on $\RR^n$, we let
this be our $\Phi$.

The coefficients of $t^j$ in the remainder are successively {\it transport
equations} for $a_0$, $a_1$, and so on. Consider the transport equation for
$a_0$. Fix a $w$ and choose normal coordinates for $z$ centred at
$z=w$. Then
$$
\Delta = \sum_j D_{z_j}^2 + O(z) D_z \text{ and } \Phi(z,w) = \abs{z}^2/2 + O(\abs{z}^3),
$$
so 
$$
\Delta_z \Phi(z,w) + n = O(z) \text{ and } 
g(\nabla \Phi, \nabla) = \sum_j (z_j + O(\abs{z}^2)) \pa_{z_j} .
$$
Thus, the transport equation for $a_0$ has the
form
$$
(z_i + O(\abs{z}^2)) \frac{\pa}{\pa z_i} a_0 = f\cdot a_0, \quad f = \half
\Delta \Phi + \frac n2
 = O(z),
$$
where all terms are smooth. This has a unique smooth solution satisfying \begin{equation}
a_0(z,z) = 1 \text{ for all } z \in X^\circ .
\label{transport-init}\end{equation}

By the stationary phase lemma (see for example \cite{Ho}, Theorem 7.7.5), if $\Phi$ and $a_0$ satisfy \eqref{eik} and
\eqref{transport-init} and all $a_j$ are supported within half the
injectivity radius of the diagonal, the initial condition \eqref{init-cond}
is satisfied.

The other transport equations take the form
$$
(z_i + O(\abs{z}^2)) \frac{\pa}{\pa z_i} a_j + j a_j = f \cdot a_j - \frac{i}{2} \Delta_z a_{j-1} - i V a_{j-1},
\quad f \text{ as above.}
$$
We inductively suppose that $a_0, \dots, a_{j-1}$ are smooth. Then there is a unique smooth
solution $a_j$ of this equation, establishing the inductive hypothesis for $a_j$. 

To define $\cU_1$, we take our formal solution $\tilde \cU_1$, and multiply
by a smooth function $\chi(z,w)$ which is equal to $1$ when $d(z,w) \leq
\iota(X,g)/4$, and equal to $0$ when $d(z,w) \geq \iota(X,g)/2$. Then we
take $\cU_1$ to be an asymptotic sum of the formal series so obtained. The
error term we decompose into $\cE_1 + \cR_1$, where $\cR_1$ is given by all
terms containing a derivative of $\chi$, and $\cE_1$ is the
remainder. Since $\tilde \cU_1$ is a formal solution, $\cE_1$ is
$O(t^\infty)$ as $t \to 0$, while $\cR_1$ is supported away from the
diagonal. This completes Step 1.

\

\subsection{Step 2 --- at the main face} 
Here we solve away the error term $\cR_1$. We regard $\cR_1$ as a half-density by multiplication by the Riemannian half-density $\abs{dg_z dg_w dt}^{1/2}$. Then $\cR_1$ has an asymptotic expansion of the form 
\begin{equation}
t^{-n/2 + 1} e^{i\Phi(z,w)/t} \sum_{j=0}^\infty t^j r_j(z, w) \abs{dg_z dg_w dt}^{1/2}.
\end{equation}
It is therefore a Legendrian distribution, associated to the Legendrian $L$ of Lemma~\ref{Leg} --- see \cite{MR96k:58230}, or section~\ref{subsec:Leg} for a brief description. 

We now need to specify how our differential operator $\so$ acts on half-densities; we shall do this by specifying a flat (i.e.\ covariant constant) half-density. One natural choice would be to specify that the Riemannian half-density $\abs{dg_z dg_w dt}^{1/2}$ is flat. However, we shall make a different choice here since we want to apply the symbol calculus from \cite{MR2002i:58037}, where it is assumed that the flat half-density is a scattering density (that is, a bounded non-vanishing section of the scattering half-density bundle). Thus, we specify instead that 
\begin{equation}
\abs{\frac{dg_z dg_w dt}{t^{2n+2}} }^{1/2} \text{ is flat, \ie } \ 
\nabla \abs{\frac{dg_z dg_w dt}{t^{2n+2}} }^{1/2}  = 0.
\end{equation}
For convenience we define
\begin{equation}
\alpha = \abs{\frac{dg_z dg_w dt}{t^{2n+2}} }^{1/2}.
\label{alpha}\end{equation}
Hence, we want to solve 
\begin{equation}
\so \big( \cR_2 \alpha \big) = - \cR_1 \alpha,
\end{equation}
or equivalently
\begin{equation}
\big( \so + i(n+1)t \big) \big( \cR_2 \abs{dg_z dg_w dt}^{1/2} \big) = - \cR_1 \abs{dg_z dg_w dt}^{1/2}.
\end{equation}
The additional term $i(n+1)t$ is a subprincipal term, since it vanishes to first order at $t=0$. 

Writing $\cR_1$ in terms of $\alpha$ and using \eqref{sf-Leg} we see that 
$\cR_1 \in I^{7/4}(L)$.
The Legendrian $L$ is characteristic for $\so$ (that is, the symbol $p$ of $\so$ vanishes on $L$, by construction of $L$), so we look for
$\cU_2 \in I^{3/4}(L; \HDsc)$ whose symbol satisfies the transport equation
\eqref{transport} along $L$. 
In our situation,  $\scH_p$ is given by \eqref{Hflow-int}, while $\pa p/\pa \tau = 1$, and the subprincipal symbol $p_{\sub} = i$ according to formula (2.9) of \cite{MR2002i:58037}. 
Therefore, by \eqref{Hflow-int}, we can solve away the error term $\cR_1 \in
I^{7/4}(L)$ with a Legendrian distribution $u_0$ of order $3/4$, by
solving the ODE
\begin{equation}
-i \Big( {\mathcal{L}}_{H_p} - \frac{n}{2} \Big) \sigma^{3/4}(u_0) =
 \sigma^{7/4}(\cR_1)
\label{corr}\end{equation}
with `initial condition' that the symbol of $u_0$ vanishes in a
neighbourhood of $\{ z=w, \zeta = \eta = 0 \}$. All bicharacteristics
originate here and tend to the boundary $\{ x = 0 \}$ by the
non-trapping assumption, so there is a unique smooth solution with this
property. It follows that in the region $x \geq \epsilon$, we have
$$
\so u_0 - \cR_1 \in I^{7/4+1}(L)
$$

Inductively, assuming that $u_k \in I^{3/4}(L)$ solves 
$$
\so u_k - \cR_1 \in I^{7/4 + k + 1}(L)
$$ in the region $x \geq \epsilon$, we can find a $u_{k+1} \in I^{3/4}(L)$
which solves
$$
\so u_{k+1} - \cR_1 \in I^{7/4 + k + 2}(L)
$$
in the region $x \geq \epsilon$, by solving
\begin{equation}
-i \left( {\mathcal{L}}_{H_p} +k-\frac n2 \right) \sigma^{3/4 +
 k + 1}(v_k) = \sigma^{7/4+k+1} \big( \so u_k - \cR_1 \big)
\label{corr2}\end{equation}
and letting $u_{k+1} = u_k + v_k$. The $v_k$ can then be asymptotically
summed to yield a $u \in I^{3/4}(L)$ which solves
$$
\so u - \cR_1 \in t^N \CI( \{ x \geq \epsilon \} \times
\cG \times \halfline) \text{ for all } N.
$$
We choose a cutoff function $\chi_2 \in \CI(X)$ which is equal to $1$ when $x \geq 2\epsilon$ and $0$ when $x \leq \epsilon$. Define
$$
\cU_2 = \chi_2 u,
$$
then $\cU_2$ satisfies 
$$
\so \cU_2 = -\cR_1 + \cR_2 + \cE_2,
$$ where $\cR_2$ is a Legendre distribution of order $7/4$ supported in
$\epsilon \leq x \leq 2\epsilon$ (the $\cR_2$ error comes from derivatives
hitting the cutoff function $\chi_2$). Our new error term $\cR_2$ is now localized near the corner.

\

\subsection{Step 3 --- near the corner}\label{subsec:step3}  Recall that we expect the fundamental solution to have
quadratic oscillations at spatial infinity, i.e.\ at $x=0.$ In fact, we
expect that $e^{-i/2tx^2} U(t)$ is a somewhat simpler kernel than $U$
itself, near infinity. Therefore, we look for $\cU_3$ in the form
$$
\cU_3 \abs{dg_z dg_w dt}^{1/2} = e^{i/2tx^2} {\cUt}_3 \abs{dg_z dg_w dt}^{1/2},
$$
where $\cUt_3$ solves
$$
 e^{-i/2tx^2}  \Big( \so e^{i/2tx^2} (\cUt_3 \alpha) \Big) = 
 e^{-i/2tx^2}  \big( -\cR_2 \alpha \big).
$$
Here $\alpha$ is the half-density \eqref{alpha}. Then in terms of 
$\cUt_3$ we have
\begin{equation}
\Big( t^2 D_t - tx D_x + t^2 H + i(n+1)t - \frac{int}{2} \Big)  \cUt_3 = - e^{-i/2tx^2} \cR_2.
\label{mod-eqn}\end{equation}

Here the term $i(n+1)t$ comes from the fact that we have chosen the scattering half-density $\alpha$ rather than the Riemannian half-density to be covariantly constant, as in Step 2. Note that the differential operators appearing in \eqref{mod-eqn} are generated by the fibred-scattering vector fields on $\XXt$ appearing in \eqref{sf-vfs}. This is not too surprising, since the analytic operation of multiplying $\cU_3$ by $e^{-i/2tx^2}$ corresponds to the geometric operation of \eqref{add} which maps the Lagrangian $\Lt$ to a smooth Legendrian $\Lphi$ in $\TstarphiXXt$. Since $\Lphi$ is well-behaved in $\TstarphiXXt$, in the sense of being smooth up to the boundary, we can expect the corresponding operator \eqref{mod-eqn} to be well-behaved in terms of the fibred-scatteringfibred-scattering structure on $\TstarphiXXt$. 

\begin{lemma} 
The kernel $e^{-i/2tx^2} \cR_2$ is an element of
$I^{7/4, \infty}(\Lphi)$, as defined in \cite{MR2002i:58037}.
\end{lemma}

\begin{proof} If $\Psi$ parametrizes the Legendrian $L$, then $\Psi -
  1/(2tx^2)$ parametrizes $\Lphi$ according to Lemma~\ref{lemma:fibredLeg}.
Hence $e^{-i/2tx^2} \cR_2$ is Legendrian with respect to $\Lphi$, of order $7/4$ at mf according to the calculation in Step 2. Since $\cR_2$ vanishes in a neighbourhood of $x = 0$, it is order $\infty$ at sf. 
\end{proof}

Thus, it makes sense to seek the solution $\cUt_3$ to \eqref{mod-eqn}
in the space $I^{3/4, *}(\Lphi)$. Here the value of $*$, specifying
the rate of decay of the symbol as $x \to 0$, is determined by the
transport equation which is regular singular at $x=0$.

To do this, we must analyze the Hamilton vector field of our operator
on $\Lphi$ near the boundary at $x=0$. This is given (modulo some
identifications) by the vector field $V$ in \eqref{hamflow:blownup},
which is $x \pa_x + x W$, where $W$ is tangent to the boundary. To
determine the subprincipal symbol, observe that the adjoint of $t(xD_x
- in/2)$ with respect to the Riemannian density is $t(xD_x - in/2) + O(tx)$ near the corner, hence the Weyl
symbol of $t(xD_x - in/2)$ is $\nu + O(tx)$ where $\nu$ is the variable dual to
$x$. Hence the subprincipal symbol of $t(xD_x - in/2)$ at mf vanishes at
the boundary of mf (that is, at $x=0$).  Also, the subprincipal symbol of $t^2 H$  vanishes identically
on the main face.  So the boundary subprincipal symbol of the operator in \eqref{mod-eqn} is equal to $i + O(x)$.  Thus, we look
for a $u_0 \in I^{3/4, *}(\Lphi)$ satisfying
\begin{equation}
\Big( -i {\mathcal{L}}_{-x \pa_x + x W} - i ( \frac1{2} + \frac 34 -
\frac{2n+1}{4} ) +i  + xq' \Big) \sigma^{3/4}(u_0) =
\sigma^{7/4}(e^{-i/2tx^2} \cR_2)
\end{equation}
This gives an equation of the form
$$ \Big( \mathcal{L}_{x \pa_x + xW} + \frac n2 + x q' \Big) \sigma^{3/4}(u_0) =
-i\sigma^{7/4}(e^{-i/2tx^2} \cR_2),
$$ so by \eqref{legendriansymbol}, $\sigma^{3/4}(u_0)$ is of the form
$a_0 x^{-n/2} \abs{x^{-1} dx d\lambda}^{1/2} \abs{dt}^{3/4 - (2n+1)/4}$,
where $a_0$ is smooth, and $\lambda$ are coordinates on the boundary
of $\Lphi$, extended into the interior. Here we have explicitly
included the power of $\abs{dt}$ in the formula since we now need to
change boundary defining function from $t$ to $\rho = tx$, in order to
apply formulae from \cite{MR2002i:58037}. Doing this, we find that
$\sigma^{3/4}(u_0)$ is of the form $a_0 x^{-1/2} \abs{x^{-1}dx
  d\lambda}^{1/2} \abs{d\rho}^{3/4 - (2n+1)/4}$, where $a_0$ is
smooth. By the (corrected) symbol calculus for fibred Legendrians given by Proposition 3.4 of \cite{MR2002i:58037}\footnote{The power of $\rho$ in the final nonzero term of the exact sequence of Proposition 3.4 of \cite{MR2002i:58037} is incorrect. It should be $\textbf{r} - m$, not $m - \textbf{r}$.}, this implies that $u_0 \in I^{3/4, 1/4}(\Lphi)$, and that 
\begin{equation}
\Big( t^2 D_t - tx D_x + t^2 H + i(n+1)t - \frac{int}{2} \Big)  u_0 +  e^{-i/2tx^2} \cR_2 \in I^{7/4 + 1, 1/4}(\Lphi).
\label{step3-1}\end{equation}
Inductively, we look for $u_l \in I^{3/4, 1/4}(\Lphi)$ such that 
\begin{equation}
\Big( t^2 D_t - tx D_x + t^2 H + i(n+1)t - \frac{int}{2} \Big)  u_l +  e^{-i/2tx^2} \cR_2 \in I^{7/4 + l+1, 1/4}(\Lphi).
\label{step3-2}\end{equation}
In fact, we will show more --- we will show that there is $u_l$ as above whose (full) symbol is of the form 
\begin{equation} a(\rho, x) x^{-1/2} \abs{x^{-1}dx
  d\lambda}^{1/2} \abs{d\rho}^{3/4 - (2n+1)/4}, \quad \text{ with $a$ smooth.}
\label{rho-reg}\end{equation}
The important point here is that $a$ is a smooth function of $\rho$ and $x$, not just a smooth function of $t$ and $x$. To show this, we use the boundary defining function $\rho$ rather than $t$ for the main face, even though it is degenerate at the corner $t=x=0$. In terms of the coordinates $\rho, x, y, w$ near the corner, partial derivatives transform as
\begin{equation}\begin{aligned}
t D_t \big|_{x,y} &= \rho D_\rho \big|_{x,y} \\
 x D_y \big|_{x,t} &= x D_y \big|_{\rho, x} \\
x D_x \big|_{y,t} &= x D_x \big|_{y, \rho} + \rho D_\rho \big|_{x,y}.
\end{aligned}\end{equation}
We also write $V = x^2 \tilde V$ where $\tilde V \in \CI(X)$ by assumption. 
Hence the operator in \eqref{step3-2} takes the form 
\begin{equation}
\tilde P = -\rho D_x + \frac{\rho^2}{2} \Big( (x D_x + \rho D_{\rho})^2 + i(n-2) (x D_x + \rho D_{\rho}) + \Delta_{h(x)} + \tilde V \Big) + i \rho \frac{n+2}{2x}.
\label{P3}\end{equation}
Any Legendrian distribution $u$ in  $I^{3/4, 1/4}(\Lphi)$ can be written, modulo $x^\infty t^\infty C^\infty(X \times \halfline)$, as a finite sum of oscillatory integrals of the form (written as a Taylor series in $\rho$)
\begin{equation}
u = \int_{K} e^{i\psi(x,y,w,v)/\rho} \rho^{n/2 + 1 - k/2} \sum_{j=0}^\infty \rho^j a_j(x,y,w,v) \, dv
\label{u3}\end{equation}
where $\psi$ parametrizes $\Lphi$, $K \subset \RR^k$ is bounded, and $a_j \in x^{-n/2 - 1 - j}C^\infty$. We claim that in fact, $a_j \in x^{-n/2 - 1}C^\infty$, or in other words, that the full symbol of $u$ is smooth in $x$ and $\rho$. 

The equation to be solved is 
\begin{equation}
\tilde Pu = -e^{-i/2tx^2} \cR_2 = \int_{K} e^{i\psi(x,y,w,v)/\rho} \rho^{n/2 + 2 - k/2} \sum_{j=0}^\infty \rho^j b_j(x,y,w,v) \, dv
\end{equation}
where the $b_j$ are supported in $x \geq \epsilon$. Substituting \eqref{P3} for $\tilde P$, we find that the left hand side is given by 
\begin{equation}\begin{gathered}
\int_K e^{i\psi(x,y,w,v)/\rho} \rho^{n/2 + 1 - k/2} \Bigg\{ \Big( -\psi_x + \half (x \psi_x - \psi)^2 + \half |d_y \psi|^2_{h(x)} \Big)  
+ \rho \Big( -D_x + \\ (x\psi_x - \psi)(x D_x + \rho D_\rho) + \langle \nabla_y \psi, D_y \rangle_{h(x)} + i \frac{n+2}{2x} + f \Big) + \rho^2 Q \Bigg\}
\sum_{j=0}^\infty \rho^j a_j(x,y,w,v) \, dv ,
\end{gathered}\end{equation}
where $Q$ is a second order operator generated over $\CI$ by $x D_x$, $\rho D_{\rho}$, and $D_y$. 
Since $\psi$ parametrizes a Legendrian which is characteristic for $\tilde P$, the expression 
\begin{equation} -\psi_x + \half (x \psi_x - \psi)^2 + \half |d_y \psi|^2_{h(x)}
\label{eikonal-3}\end{equation}
vanishes when $d_v \psi = 0$. As $\psi$ satisfies the nondegeneracy condition \eqref{ident-fib}, the function \eqref{eikonal-3} can be expressed
$$
\sum_{l=1}^k a_l \frac{\pa \psi}{\pa v_l} \quad \text{ with $a_l$ smooth.}
$$
We then write $a_l (\pa_{v_i} \psi ) e^{i\psi/\rho} = a_l \rho D_{v_i} e^{i\psi/\rho}$ and integrate by parts in $v$. These terms then become $O(\rho)$ terms.  Let us write $W = (x\psi_x - \psi)x D_x + \langle \nabla_y \psi, D_y \rangle_{h(x)} + \sum_l a_l D_{v_l}$.  The equation becomes  
\begin{equation}
\sum_{j=0}^\infty \rho^j 
\Bigg\{ \Big( -D_x + W + i \frac{n+2}{2x} + f_j  \Big) + \rho Q_j \Bigg\} a_j 
= \sum_{j=0}^\infty \rho^j b_j. 
\end{equation}
Here, $f$ and $f_j$ are smooth and $Q_j$ is a second order operator generated by $x D_x, D_y$ and $D_v$ with smooth coefficients. Thus, we need to solve
\begin{equation}\begin{aligned}
(-xD_x + xW +  i \frac{n+2}{2} + xf_0)a_0 &= xb_0 , \\
(-xD_x +xW +  i \frac{n+2}{2} + xf_j)a_j &= xb_j + xQ_j a_{j-1}.
\end{aligned}\label{trans}\end{equation}
Here $b_j$ is supported in $x \in [\epsilon, 2\epsilon]$ and our initial condition for $a_j$ is that it vanish for $x > 2\epsilon$. These are regular singular ODEs for the $a_j$. 
The first equation certainly has a solution which is in $x^{-n/2-1} C^\infty$. Assuming that $a_j$ is in $x^{-n/2-1} C^\infty$, it follows from \eqref{trans} that $a_{j+1} \in x^{-n/2-1} C^\infty$, since the right hand side of \eqref{trans} is in $x^{-n/2} C^\infty$ and the indicial root of the regular singular operator on the left hand side is $-n/2 - 1$. Hence the claim follows by induction. (Notice that this induction would not work if we did not make the stronger assumption \eqref{rho-reg} on the symbol of $u$.)

Thus, we let $\tilde \cU_3$ be an asymptotic sum of the formal series \eqref{u3} constructed above, supported say in $t \leq 1$. This is certainly in $I^{3/4, 1/4}(\Lphi)$ and solves the equation
\begin{equation}
\Big( t^2 D_t - tx D_x + t^2 H + i(n+1)t - \frac{int}{2} \Big)  \cUt_3 = - e^{-i/2tx^2} \cR_2 + \cE_3 +  \tilde \cR_3,
\label{error-3}\end{equation}
where $\tilde \cR_3 \in I^{\infty, 1/4}(\Lphi)$. This completes Step 3.

\


\subsection{Step 4 --- near spatial infinity}\label{subsec:step4} Now we have an error term $\tilde \cR_3 \in 
I^{\infty, 1/4}(\Lphi)$, and we seek a solution to
\begin{equation}
\Big( t^2 D_t - tx D_x + t^2 H + it \frac{n+2}{2} \Big)  \cUt_4 = -  \tilde \cR_3 \in x^\infty t^\infty \CI(\XXt).
\label{eqn-4}\end{equation}
We seek a solution $\cUt_4$ in the space $I^{\infty, 1/4}(\Lphi)$.

We can write $-\cRt_3$ as a finite sum of terms of the form 
\begin{equation}
\int_{K} e^{i\psi(x,y,w,v)/\rho} \rho^{-k/2} \sum_{j=0}^\infty \rho^j b_j(t,y,w,v) \, dv
\label{u4}\end{equation}
where we have expanded the symbol as a Taylor series in $\rho$. Here each $b_j$ is $O(t^\infty)$ at $t=0$. In terms of the coordinates $(t, \rho, y),$
$$
t D_t \big|_{\rho,y} = t D_t \big|_{x,y} - x D_x \big|_{t,y},
$$
so the operator in \eqref{eqn-4} takes the form  
\begin{equation}
 t^2 {D_t} \big|_{\rho}  + it \frac{n+2}{2} + (\rho^2 D_\rho)^2 + h^{ij} (\rho
  D_{y_i})(\rho D_{y_j}) + \frac{\rho}{t} Q 
\label{cUt4-2}\end{equation}
where $Q$ is a second order differential operator generated over $C^\infty$ by the
vector fields $\rho^2 D_\rho$ and $\rho D_y$. Hence, we wish to solve 
\begin{equation}\begin{gathered}
\!\!\! \Big( t^2 {D_t} \big|_{\rho}  + it \frac{n+2}{2} + (\rho^2 D_\rho)^2 + h^{ij} \rho^2
  D_{y_i} D_{y_j} + \frac{\rho Q}{t}  \Big) \!  \int\limits_{K} \! e^{i\psi(x,y,w,v)/\rho}   \sum_{j=0}^\infty \rho^{j-\frac{k}{2}} a_j(t,y,w,v) \, dv  
\\ =   \int_{K} e^{i\psi(x,y,w,v)/\rho} \sum_{j=0}^\infty \rho^{j-\frac{k}{2}} b_j(t,y,w,v) \, dv .\end{gathered}\label{series-4}\end{equation}
This gives us equations
\begin{equation}\begin{aligned}
\Big( t^2 D_t + it \frac{n+2}{2} + c(y,w,v) \Big) a_0 &= b_0, \\
\Big( t^2 D_t + it \frac{n+2}{2} + c(y,w,v) \Big) a_j &= b_j + R_j(a_0, \dots , a_{j-1}),
\end{aligned}\end{equation}
where $c(y,w,v) = \psi(0,y,w,v)^2 + h^{ij} \psi_{y_i}(0,y,w,v)\psi_{y_j}(0,y,w,v)$ is independent of $t$. The operator $R_j$ is such that $t^{j+1} R_j$ is a second order differential operator in $y$ with smooth coefficients. We show inductively that there is a solution with each $a_j \in t^\infty C^\infty$. Indeed, the equation for $a_0$ is explicitly solved by 
$$
a_0(t) = t^{(n+2)/2} e^{ic/t} \int_0^t i e^{-ic/s} b_0(s) s^{-n/2 + 3} \, ds
$$
which is in $t^\infty C^\infty$ since $b_0 \in t^\infty C^\infty$. Inductively assuming that $a_0, \dots , a_{j-1}$ are in $t^\infty C^\infty$, then it follows from the formula above, with $b_0$ replaced by $b_j + R_j(a_0, \dots , a_{j-1})$ that $a_j \in t^\infty C^\infty$. Let $\cUt_4$ be an asymptotic sum of the series on the left hand side of  \eqref{series-4}. Then $\cUt_4$ solves \eqref{eqn-4} up to an $x^\infty t^\infty C^\infty$ error term. This completes the construction of the parametrix in the nontrapping case.

\subsection{Exact solution}

We have now shown that
\begin{equation}
(D_t + H) \cU  (x, y, w, t) = e(x,y,w, t) \in \CdotI (\halfline \times
X\times \cG)= \CdotI(\halfline; \CdotI(X\times \cG)).
\label{inhomog}
\end{equation}
The exact propagator $U(t)$ is given in terms of the parametrix
$\cU(t)$ by Duhamel's formula
\begin{equation}
U(t) = \cU(t) + i\int_0^t U(s)\big( D_t + H) \cU(t-s) \, ds.
\label{dH}\end{equation}

It follows from a commutator argument due to Craig, \cite{MR99f:35028b}
Th\'eor\`eme~14, that, letting
$$
\cH_k = \bigcap_{s=0}^k x^{k-s} \Hsc[s](X),
$$
we have for all $k>0$
\begin{equation}
U(t): \cH_k \to L^\infty_{\loc}(\RR_t; \cH_k).
\label{prop-map}\end{equation}
Note that $\bigcap_k \cH_k = \CdotI$, so in particular
\begin{equation}
U(t) : \CdotI(X) \to \CdotI(X) \ \text{ for all } t \leq 0.
\label{Schw-map}\end{equation}

\begin{lemma}\label{schw}
 Let $e \in \CdotI(\halfline; \CdotI(X \times \mathcal{G}))$. Then    $$
Ke(t)(z,w) \equiv \int_0^t  ds \int_X U(s)(z,z') e(t-s)(z',w) \, dg(z')   \in \CdotI(\halfline; \CdotI(X \times \mathcal{G})).
$$
\end{lemma}

\begin{proof} 
Let us fix $w \in \mathcal{G}$; then we may regard $e(z, w, t)$ as an element of $t^n \CdotI(\halfline; \cH_k)$ for any $n, k \in \NN$. 
First we establish that $Ke \in t^n L^\infty(\halfline; \cH_k)$. We compute for $t \in [0, T]$
\begin{align*}
\norm{Ke(t)}_{\cH_k} &\leq \int_0^t \norm{U(s) e(t-s)}_{\cH_k} \, ds
\\
&\leq C_T \int_0^t \norm{e(t-s)}_{\cH_k}\, ds\\
&=O(t^n) \text{ for all } n\in \NN.
\end{align*}
This shows that  $Ke \in t^n
L^\infty(t; \cH_k)$ for every $k$ and $n$. Derivatives of $t$ and $w$ can now be estimated similarly and uniform estimates then follow from the relative compactness of $\mathcal{G}$. 
\end{proof}

We conclude from \eqref{dH} and Lemma~\ref{schw} that $U(t) - \cU(t)
\in \CdotI(\halfline; \CdotI(X \times \mathcal{G}))$. This proves that
the kernel of $U(t)$ is such that $e^{-i/2x^2t}U(t) \in I^{3/4,
  1/4}(\Lphi)$, which completes the proof of Theorem~\ref{main}.


\section{The trapping case}\label{section:trapping}

We now sketch briefly the changes to be made in the parametrix construction
in case there exist trapped geodesics.  In this case, given a properly
supported $Z \in \Psi^0 (\cG)$ with $\WF' Z \subset \nft,$ the
non-forward-trapped part of the phase space, we wish to
construct a partial parametrix $\cUZ$ such that
$$
(D_t + H) \cUZ  (z,w, t) = e(z,w, t) \in \CdotI(\halfline; \CdotI (X\times \mathcal{G})).
$$
with initial condition
$$
\cUZ(t) f \to Zf \text{ in } C^{-\infty}(X) \text{ as } t \to 0.
$$
It then follows by Duhamel's principle that $\cUZ - U(t)Z
\in \CdotI(\halfline; \CdotI (X\times \mathcal{G})).$

The only modification necessary in the construction comes at the end of Step 1. We begin with a result about the
composition of pseudodifferential operators in the boundary with
the simple Legendrian distributions appearing in Step 1.

\begin{lemma}
Let $k(z,w,t)=t^s a(z,w,t) e^{i \Phi(z,w)/t}$ such that $a$ and $\Phi$ are smooth, $a$ is compactly supported and $d_w \Phi \neq 0$ on the support of $a$.  Let $B \in
\Psi^0(X)$ be a properly supported, classical pseudodifferential operator.  Then
$$
v=\int_{X} k(z,w,t) B(w, w') \, dg(w)
$$
has  the form
$$
v = t^s \tilde a(z, w, t) e^{i\Phi(z,w)/t}, \quad \tilde a \text{ smooth, }
$$
with 
$$
\WFsc(v) \subset \WFsc(k) \cap \{ (z,w, \zeta, \eta, \tau) \mid (w, \eta) \in \WF' (B) \}.
$$
\label{lemma:composition}
\end{lemma}
\begin{proof}
The result follows directly from the lemma of stationary phase. Alternatively, we may regard $z$ and $t$ as smooth parameters and regard $v$ as $B^T (k)$, the transpose of $B$ acting on the distribution $k$. Recall that
when $B$ is a zeroth order classical pseudodifferential operator and $c,
\phi \in \CI$ with $c$ compactly supported and $d\phi \neq 0$ on $\supp c,$
$$
B \big( c(w) e^{i\phi(w)/t} \big)  = \tilde c(w,t) e^{i\phi(w)/t} 
$$
where $\tilde c\in \CI(\halfline \times X)$ and $\tilde c(w,t) =
O(t^\infty)$ for all $w$ such that $(w, d_w\phi) \notin
\WF'B.$  (See \cite{MR82i:35172}, Chapter 8, section 7.) The lemma now follows directly. 
\end{proof}

To construct $\cUZ$, we construct $\cU_1$ as in Step 1 of the previous section and consider $\cU_1 \circ Z$. The error term $\cR_1 \circ Z = (D_t + H) \cU_1 \circ Z$ is then a Legendre distribution associated to $L$, by the lemma just proved (since $d_w \Phi \neq 0$ on the support of the symbol of $\cR_1$). Moreover the lemma shows that the scattering wavefront set of 
$(D_t + H) \cU_1 \circ Z$ lies on the portion of $L$ emanating from $(w, \eta) \in \WF'(Z)$, which is contained in the non-forward trapped set. Now we proceed with Steps~2, 3 and 4 with
$\cR_1 \circ Z$ substituted for $\cR_1$.  As no trapped rays lie in the
microsupport of $\cR_1 \circ Z$ however, we need only solve the transport
equation in Step 2 along the remaining, non-trapped, rays.  The rest of the
construction proceeds as before.

Thus, we have shown that $\cUZ$ is the sum of two terms, one of which has the form $\cU_1 \circ Z$ where $\cU_1$ is a Legendrian distribution as in subsection~\ref{Step1}, and the other is a fibred-scattering Legendrian distribution of order $(\frac{3}{4}, \frac{1}{4})$. 

{\it Remark.}  It is not true that $\cU_1\circ Z$ is a Legendrian
distribution associated to $L$, since $d_z \Phi = 0$ at $z=w$. In fact,
$\cU_1\circ Z$ is a Legendrian conic pair associated to $(\Delta, L)$, where
$\Delta$ is the diagonal Legendrian $\{ (z, \zeta, z, -\zeta, \kappa = 0)
\}$. Since we do not need this fact, we omit the proof.


\section{Egorov theorem}
Let $W(t) = e^{-i/2tx^2} U(t).$ Suppose $A$ is a properly supported pseudodifferential
operator on $X^{\circ}$, microsupported in $O \times O$ with $O \subset \nbt$.  We will 
show that $B = W(t) A W(t)^*$ is a scattering pseudodifferential operator
of the same order as $A$, and that the symbol and microsupport of $B$ are
given by the pullback of the symbol and microsupport of $A$ by the sojourn
relation, scaled by the factor $t^{-1}$.  Indeed we prove a more general result. To state it, suppose that $S$ is a contact diffeomorphism from $\Openset \subset S^* X^{\circ} \to \Tscstar[\pa X] X$. Then, by definition, $S^* \chi$, the pullback of the contact form on $\Tscstar[\pa X] X$, is equal to a multiple $f \hat \eta \cdot dw$ of the contact form on $S^* X$. Associated to $S$ is a Legendrian submanifold $G$ of $\Tscstar[\pa X \times X^\circ] X \times X^\circ$, given by 
\begin{equation}
G = \{ (y, w, \nu, \mu, -f \eta) \mid (w, \eta) \in \Openset, (y, \nu, \mu) = S(w, \eta) \}.
\label{graph}\end{equation}
Indeed this is the correspondence of Lemma~\ref{lemma:contactomorphism}. 
\begin{proposition}\label{Egorov}
Let $S$ be a contact diffeomorphism $S : \Openset \subset S^* X^{\circ} \to \Tscstar[\pa X] X$ and let $G$ be the Legendre submanifold given by \eqref{graph}. Suppose that 
$W\in I^0(X \times X^\circ ,G; \Omegasc^{\h} \times \Omega^{\h})$, and suppose that $A\in \Psi^m (X^\circ; \Omega^{\h})$ is properly supported and microsupported in $\Openset$. Then
 \begin{equation}
  W A W^* \in \Psisc^{-\infty, m} (X; \Omegasc^{\h}) \text{ with symbol }
 \abs{\sigma(W)(q_G)}^2  \cdot \sigma (A)(S^{-1}(q))
\label{WAWstar}\end{equation}
where $q_G$ is the point on $G$ corresponding to $q$. 
Conversely, if $\tilde A \in \Psisc^{-\infty, m}(X; \Omegasc^{\h})$ has compact microsupport contained in the range of $S$, then modulo a kernel in $\CdotI(X \times X)$, $W^* \tilde A W \in \Psi^m(X^\circ)$ is a properly supported pseudodifferential operator with symbol 
\begin{equation}
\sigma(W^* \tilde A W)(\zeta) = \sigma(\tilde A)(S(\zeta)) \cdot \abs{ \sigma(W)(\zeta_G)}^2.
\label{WstarAW}\end{equation}
\end{proposition}

The proof proceeds as follows: we reduce to the case $X=\overline{\RR^n},$
the radial compactification of Euclidean space, by localization; applying
the Fourier transform and results of Melrose-Zworski \cite{MR96k:58230}, we
then deduce the result from the standard Egorov theorem.

To begin, we need a lemma, in effect a version of Proposition 10 of
\cite{MR96k:58230} with parameters, which tells us what happens when we
pull back a Legendrian distribution on $\RR^n \times \overline{\RR^n}$ by
the Fourier transform $\Fourier.$  The Fourier transform invariantly
maps half-densities on a vector space to half-densities on the dual space. We recall from \cite{MR96k:58230} that it can be interpreted as a Fourier integral operator associated to the Legendre diffeomorphism $L : \Tscstar_{\pa \overline{\RR^n}} \overline{\RR^n} \to S^*((\RR^n)^*)$ given by 
$L(\hat z, \zeta) = (\zeta, -\hat z)$, where $\hat z \in S^{n-1}$, the radial compactification of $\RR^n$, and $\zeta \in (\RR^n)^*$, the dual space. 

\begin{lemma}
  Let $S$, $G$ and $W$ be as in Proposition~\ref{Egorov} for $X = \overline{\RR^n}$, and let $\Fourier$ be the Fourier transform.   Then
  $$
 \Fourier \circ W \in I^{0}((\RR^n)^* \times \RR^n, G_L;
  \Omega^{\h}) + \mathcal{S}((\RR^n)^* \times \RR^n),
  $$
where $G_L$ is the Lagrangian associated to the graph of the contact transformation $L \circ S$, $L$ is the Legendre diffeomorphism, and $\mathcal{S}$ denotes Schwartz space.
\label{lemma:FT}
\end{lemma}

\begin{proof}[Proof of Lemma~\ref{lemma:FT}]
Without loss of generality we may assume that $W$ can be written in terms of a single parametrization. Then $W$ can be written 
$$
(2\pi)^{-\frac{k}{2}} \int_{\RR^k} e^{i{\phi(\hat z,w,v)}|z|}
    a(\hat z, \frac1{|z|}, w ,v) |z|^{\frac{k}{2}} \, dv \, \abs{dz \, dw}^{\h}, \quad a \text{ smooth, compactly supported,}
$$
modulo $\mathcal{S}(\RR^n \times \RR^n)$, so
$\Fourier \circ W$ is given by
$$
(2 \pi)^{-n/2-k/2} \int \!\! \int e^{-i z \cdot \zeta
  +i{\phi(\hat z,w,v)}|z|}  a(\hat z,|z|^{-1},w,v) |z|^{k/2}  \, dv \, dz \,
  \abs{d\zeta}^{\h}\abs{dw}^{\h}
$$
modulo $\mathcal{S}((\RR^n)^* \times \RR^n)$.
Imitating the proof of Proposition 10 of \cite{MR96k:58230}, we set $\tilde v=v |z|$ so that $\Fourier \circ W $ is given, modulo $\mathcal{S}((\RR^n)^* \times \RR^n)$, by
\begin{equation}
 (2 \pi)^{-n/2-k/2} \int \!\! \int e^{-i z \cdot \zeta
  +i\phi(\hat z, w, \tilde v/|z|)|z|} a(\hat z, |z|^{-1}, w, \tilde v/|z|) |z|^{-k/2}
\, d\tilde v
  \, dz\,  \abs{d\zeta}^{\h}\abs{dw}^{\h}.
\label{lagr}\end{equation}
The phase is now homogeneous of degree one in $(z, \tilde v)$ and parametrizes the Lagrangian $G_L$ while the amplitude $a(\hat z, |z|^{-1}, w, \tilde v/|z|) |z|^{-k/2}$ has
order $-k/2$ with $k+n$ phase variables $z, \tilde v$. Hence \eqref{lagr} is a Lagrangian distribution of  
Lagrangian order $0$.
\end{proof}

\begin{proof}[Proof of Proposition~\ref{Egorov}
]
By a microlocal partition of unity, we may assume that $\opWF A \subset O$, a
small subset of $S^* X^\circ$.  Then $S(O)$ lies over a small
closed set in $\pa X$. Without loss of generality, we may assume
that the support of the kernel $W(z,w)$ in the $z$ variable lies within a small
neighborhood $W$, covered by a single coordinate chart $(x,y)$, of this closed set, since introducing a spatial cutoff introduces only
residual terms in $W A W^*$. By identifying $W$ with a neighbourhood of a point of the boundary of $\overline{\RR^n}$ and identifying a neighbourhood of $\pi(O)$ with a neighbourhood of the origin in $\RR^n$, we may assume that that $A \in \Psi^0(\RR^n)$ is properly supported and that $W$ is as in Lemma~\ref{lemma:FT}. 

Now we write
$$
WAW^* = \Fourier^* (\Fourier W) A (\Fourier W)^* \Fourier,
$$
and apply the standard Egorov theorem to the middle three factors to conclude that this takes the form
$$
\Fourier^* \tilde A \Fourier,
$$
where $\tilde A \in \Psi^m(\RR^n) + \mathcal{S}(\RR^n \times \RR^n)$ is, up to a Schwartz kernel, a pseudodifferential operator with principal symbol
$$
\tilde a(q) = \sigma(A)(L^{-1}(S^{-1}(q))) \cdot  \sigma((\Fourier W)(\Fourier W)^*)(q) = \sigma(A)(L^{-1}(S^{-1}(q))) \cdot  |\sigma(W)(q_G)|^2.
$$ 
Proposition~8 of \cite{MR96k:58230} now shows that
$$
\Fourier^* \tilde A \Fourier \in \Psisc^{-\infty, m} (\overline{\RR^n})
$$ 
with principal symbol
$$
\sigma(A)(S^{-1}(q)) \cdot |\sigma(W)(q_G)|^2.
$$
This proves \eqref{WAWstar}. The proof of \eqref{WstarAW} is similar, proceeding from the expression
$$
W^* \tilde A W = W^* \Fourier^* \Fourier \tilde A \Fourier^* \Fourier W.
$$
\end{proof}

\section{Wavefront set bound}

In this section we prove Theorem~\ref{WF-thm}. First we prove a preliminary result on the propagator at a fixed time $t \neq 0$. We define the contact transformations $\contact_t$ for $t \neq 0$ as follows: for $t>0$ and $q \in S^* X^\circ,$ let $\contact_t(q)$ with domain $\nft$ be given by $\contact_t(q)=t^{-1} S_f(q)$ where
$S_f$ is the contact transformation defined in
Lemma~\ref{lemma:contactomorphism} and the scaling by $t^{-1}$ acts in the
fibre variables.  For $t < 0$, let $\contact_{t}$ with domain $\nbt$ be given by $\contact_{t}(q)=|t|^{-1}S_b(q)$. Let $G_t$, $t \neq 0$, be the Legendrian \eqref{graph} determined by $S_t$. 

\begin{lemma}\label{Wt} 
Let $Z \in \Psi^0(X^{\circ})$ be properly supported, and microsupported inside $\nft$. 
For fixed $t > 0$, the operator $W(t)Z$ is in $I^0(X \times X^{\circ}; G_t)$. Similarly, if $Z' \in \Psi^0(X^{\circ})$ is properly supported, and microsupported inside $\nbt$, then for fixed $t < 0$, the operator $W(t)Z'$ is in $I^0(X \times X^{\circ}; G_t)$.
\end{lemma}

\begin{proof} The construction of the parametrix shows that $W(t) Z$ is given by a finite sum of oscillatory integrals
$$
\int\limits_{K \Subset \RR^k} e^{i\psi(0,y,w,v)/xt} x^{-k/2} a(x,y,w,v,t) \, dv \big| \frac{dx \, dy \, dw}{x^{n+1}} \big|^{1/2},
$$
where $\psi$ and $a$ are constructed in sections~\ref{subsec:step3} and \ref{subsec:step4}. (Here, because $t > 0$, we may replace $\psi(x, y, w, v)$ by $\psi(0, y, w, v)$ in the phase.) According to the remark below equation \eqref{fibred-term}, this is a Legendre distribution associated to the Legendre submanifold $G_t$, of order zero. The proof for $t < 0$ is similar. 
\end{proof}

The proof of Theorem~\ref{WF-thm} is now straightforward.

\begin{proof}[Proof of Theorem~\ref{WF-thm}]  We only prove the theorem for $t>0$, since the argument is similar for $t<0$. Let $W(t) = e^{-i/2tx^2} U(t)$ as
above. Then since both $U(t)$ and the multiplication operator $e^{-i/2tx^2}$ are
unitary and map $\CdotI(X)$ to itself, the same is true of $W(t)$. For
simplicity we shall first assume that the metric is nontrapping.

Suppose that $q \notin \WF(u(\cdot, t))$ for some $t > 0$. Then there is a properly supported $A \in \Psi^0(X)$ which is elliptic at $q$ and such that $A u(\cdot, t) \in \CdotI(X)$. 

Since $f = U(-t) u$, we have $f = e^{-i/2tx^2} W(-t)u$, or $u = W(-t)^* e^{i/2tx^2} f$. Hence 
$$
A u = A W(-t)^* e^{i/2tx^2} f \in \CdotI(X).
$$ Applying $W(-t)$, we see that also $W(-t) A W(-t)^* e^{i/2tx^2} f \in
\CdotI(X)$.  By Proposition~\ref{Egorov}, $W(-t) A W(-t)^* = \tilde A$ is a
scattering pseudodifferential operator of order zero, which is elliptic at
$\contact_{-t}(q)$. We have just shown that $\tilde A (e^{i/2tx^2} f) \in
\CdotI(X)$, so by definition of the scattering wavefront set, we see that
$\contact_{-t}(q) \notin \WFsc(e^{i/2tx^2} f)$.  Conversely, assume that
$\contact_{-t}(q) \notin \WFsc(e^{i/2tx^2} f)$. 
Then there exists a
scattering pseudodifferential operator $\tilde A$ with compact microsupport which is elliptic at
$\contact_{-t}(q)$ and such that $\tilde A (e^{i/2tx^2} f) \in
\CdotI(X)$.  By Proposition~\ref{Egorov}, we may write $\tilde A = W(-t) A W(-t)^*$, with $A = A' + A'' \in \Psi^0(X^\circ) + \CdotI(X \times X)$, where $A'$ is properly supported and elliptic at $q$ and $A''$ has a Schwartz kernel.  Hence $W(-t) A W(-t)^* e^{i/2tx^2} f = W(-t) A
u(\cdot, t)$ is in $\CdotI(X)$, which implies (using \eqref{Schw-map}) that $A u(\cdot, t) \in
\CdotI(X)$. Hence $q \notin \WF(u(\cdot, t))$.

Now we do not assume that the metric is nontrapping, but assume instead
that the point $q$ is not backward-trapped. Then since the set $\nbt$ of
non-backward-trapped points is open, one can choose $A$ as above so that
$\opWF A\subset \nbt$. In fact, one can find in addition a properly
supported pseudodifferential operator $Z$ such that $Z$ is equal to the identity microlocally on $\WF'(A)$ and such that $\opWF Z\subset \nbt$. Then $ZAZ u(\cdot,
t)\in\CdotI(X)$. Repeating the argument above, we see that
$$
W(-t) Z A Z W(-t)^* e^{i/2tx^2} f \in \CdotI(X).
$$
In the trapping case, $W(-t) Z$ is a Legendre distribution associated to
$\contact_{-t}$, so by Proposition~\ref{Egorov}, $W(-t) Z A Z W(-t)^*$ is a
scattering pseudodifferential operator which is elliptic at
$\contact_{-t}(q)$, and the rest of the argument goes as before.
\end{proof}

\end{document}